# Most unexposed taut one-relator presentation 2-complexes are finitely unsplittable.

by Fredric D. Ancel and Pete Sparks

**Abstract.** The main result of this article is that among the family of one-relator presentation 2-complexes that might be expected to be *finitely unsplittable* (not the union of two proper subpolyhedra with finite first homology groups) almost all have this property. Included among these one-relator presentation 2-complexes are all *generalized dunce hats.* A *generalized dunce hat* is a 2-dimensional polyhedron created by attaching the boundary of a disk $\Delta$ to a circle J via a map $f : \partial\Delta \to J$ with the property that there is a point $v \in J$ such that $f^{-1}(\{v\})$ is a finite set containing at least 3 points and f maps each component of $\partial\Delta - f^{-1}(\{v\})$ homeomorphically onto $J - \{v\}$. The fact that generalized dunce hats are finitely unsplittable undermines a strategy for proving that the interior of the Mazur compact contractible 4-manifold M is *splittable in the sense of Gabai* (i.e., int(M) = U ∪ V where U, V and U ∩ V are each homeomorphic to Euclidean 4-space).

**1. Introduction.** An open (non-compact boundaryless) n-manifold is *splittable in the sense of Gabai* (or, more briefly, *splittable*) if it is the union of two open subsets U and V such that U, V and U ∩ V are homeomorphic to Euclidean n-space. It is easily seen that Euclidean n-space is splittable. What is not obvious is whether there exist any other splittable contractible open manifolds, and whether there exist any non-splittable contractible open manifolds.

In 2009 David Gabai surprised the geometric topology community by showing that the *Whitehead 3-manifold* is splittable [G]. The Whitehead 3-manifold is one of the fundamental examples in 3-manifold topology; it is a contractible open 3-manifold that is not homeomorphic to Euclidean 3-space. (It was discovered by J. H. C. Whitehead in 1935 as a fatal obstruction to his own proposed proof of the 3-dimensional Poincare Conjecture. [W]) Gabai's observation inspired alternative proofs and generalizations by other topologists: in [GRW] it is proved that there exist uncountably many topologically distinct splittable contractible open 3-manifolds as well as uncountably many topologically distinct non-splittable contractible open 3-manifolds. Pete Sparks in his 2014 Ph.D. thesis [S] extended the study of the spittability phenomenon to dimension 4 by constructing an uncountable family of topologically distinct splittable contractible open 4-manifolds. The question of whether there exists a non-splittable contractible open 4-manifold remains unresolved.

The Mazur 4-manifold [M] is a compact contractible 4-manifold with boundary whose interior is not homeomorphic to Euclidean 4-space. It is a fundamental object in the study of 4-dimensional manifolds. It is not known whether the interior of the Mazur 4-manifold is splittable, and resolving this question is the primary motivation behind the results in this paper.  One strategy for attacking this question relies on the fact that the Mazur 4-manifold has a *spine* that is a *dunce hat*.  In general, a subpolyhedron P of a compact piecewise linear manifold M is called a *spine* of M if M collapses to P.  The *dunce hat* is a 2-dimensional polyhedron that is contractible but not collapsible, and according to [Z], the Mazur 4-manifold has a spine which is a dunce hat.  In general, if a compact piecewise linear manifold has a spine that can be expressed as the union of two collapsible proper subpolyhedra that intersect in a collapsible subpolyhedron, then these subpolyhedra can be thickened to give a splitting of the interior of the manifold. (See Proposition 3.1.9 of [S].)  Thus, a potential strategy for proving the splittability of the interior of the Mazur 4-manifold would be to exhibit the dunce hat as the union of two collapsible subpolyhedron that intersect in a collapsible subpolyhedron.  A consequence of the theorem proved in this paper is that this strategy can't work, because the dunce hat can't be divided into two collapsible subpolyhedra.  So the interior of the Mazur 4-manifold remains a potential candidate for a non-splittable contractible open 4-manifold.

One might hope to prove that because the Mazur 4-manifold has a spine that can't be split into two collapsible subpolyhedra, then its interior is not splittable. However, there are obstacles to giving such a proof because the relationship between compact contractible manifolds and their spines is complex.  Two non-homeomorphic compact contractible 4-manifolds can have homeomorphic spines, and a single compact contractible 4-manifold can have two non-homeomorphic spines – one of which can be divided into collapsible subpolyhedra while the other can't.  Thus, for compact contractible manifolds, a spine which can't be divided into collapsible subpolyhedra does not compel a non-splittable interior.  Specifically, a 4-ball whose interior is splittable has an disc spine, which can obviously be divided into two collapsible subpolyhedra that intersect in a collapsible subpolyhedron, and it has a dunce hat spine which, according to the results of this paper, can't be divided in this fashion.  Note that the preceding observation has the obvious consequence that the 4-ball and the Mazur 4-manifold, though they are not homeomorphic, have spines which are homeomorphic – namely, dunce hats.

We remark that in dimensions n ≥ 5, there exist splittable contractible open n-manifolds besides Euclidean n-space: the interior of every compact contractible n-manifold is splittable and every Davis n-manifold is splittable [AGS].  However, for n ≥ 5, no example of a non-splittable contractible open n-manifolds is known.

The Theorem of this paper generalizes results of [GST] where it is shown that a specific triangulation of the *original dunce hat* can't be expressed as the union of two collapsible subcomplexes.  (The original dunce hat arises when the map f : $\partial\Delta \to J$ is homotopic to a homeomorphism and the set $f^{-1}(\{v\})$ contains exactly 3 points.)

During the preparation of this manuscript, the article [Bo] was posted. The results of [Bo] significantly overlap the Theorem of this paper, but neither includes the other. The techniques of [Bo] are quite different from the techniques used here. Furthermore, the motivation for the work in [Bo] (Lusternik-Schnirelmann category issues and the Andrews-Curtis conjecture) appears at first glance to be quite separate from the motivation for this paper (the splittability of the contractible open 4-manifolds), but the connection between these questions is provocative.

The authors thank Professor Craig Guilbault for helpful discussions and an idea that played a key role in the proof of Lemma 4.

## 2. The Theorem and its Corollary.

**Definition.** Suppose $\Delta$ is a 2-dimensional disk, $J = J_1 \vee J_2 \vee \ldots \vee J_n$ is a wedge of n circles with wedge point v, and $f : \partial\Delta \to J$ is a map such that $f^{-1}(\{v\})$ is a finite subset of $\partial\Delta$. Call the quotient space $\Delta \cup_f J$ a *one-relator presentation 2-complex.* Let $q : \Delta \to \Delta \cup_f J$ denote the quotient map. Call the point v the *vertex* of $\Delta \cup_f J$.

**Definition.** Suppose $\Delta \cup_f J$ is a one-relator presentation 2-complex as above. If f maps each component of $\partial\Delta - f^{-1}(\{v\})$ homeomorphically onto $J_i - \{v\}$ for some i between 1 and n, then we call $f : \partial\Delta \to J$ a *taut* map and we call $\Delta \cup_f J$ a *taut* one-relator presentation 2-complex.

**Definition.** Suppose $\Delta \cup_f J$ is a taut one-relator presentation 2-complex as above. Let $1 \leq i \leq n$. If $f^{-1}(J_i - \{v\})$ is empty, then call $J_i$ a *naked edge* of $\Delta \cup_f J$. If $f^{-1}(J_i - \{v\})$ is a single component of $\partial\Delta - f^{-1}(\{v\})$, then call $J_i$ a *free edge* of $\Delta \cup_f J$. If $J_i$ is either a naked edge or a free edge of $\Delta \cup_f J$, then call $J_i$ an *exposed edge* of $\Delta \cup_f J$. A taut one-relator presentation 2-complex that has no exposed edges is said to be *unexposed.* Thus, $\Delta \cup_f J$ is *unexposed* if and only if for every i between 1 and n, $f^{-1}(J_i - \{v\})$ is has at least two components.

**Definition.** Call a taut one-relator presentation 2-complex *finitely splittable* if it is the union of two proper subpolyhedra each of which has a finite first homology group. Otherwise, call it *finitely unsplittable.*

We now state the principal result of this paper.

**Theorem.** Suppose $\Delta \cup_f J$ is an unexposed taut one-relator presentation 2-complex where $\Delta$ is a 2-dimensional disk, $J = J_1 \vee J_2 \vee \ldots \vee J_n$ is a wedge of n circles with vertex v, and $f : \partial\Delta \to J$ is a taut map. Then $\Delta \cup_f J$ is finitely unsplittable except in the single case in which n = 1, $f^{-1}(\{v\})$ has two elements and $f : \partial\Delta \to J$ is null-homotopic. In this exceptional case, $\Delta \cup_f J$ is a 2-sphere with a pair of points identified which we denote $\mathbb{S}^2/\{2\ points\}$.

As stated in the introduction, this theorem has consequences for the attempt to prove that the interior of the Mazur 4-manifold is splittable in the sense of Gabai. In particular, the Theorem applies transparently to the dunce hat spine of the Mazur 4-manifold and to a more general class of 2-dimensional polyhedra called generalized dunce hats. We now explain this application.

**Definition.** A *generalized dunce hat* is a taut one-relator presentation 2-complex $\Delta \cup_f J$ in which $J$ is a single circle and $f^{-1}(v)$ contains at least 3 points.

**Corollary.** Every generalized dunce hat is finitely unsplittable.

**Proof.** Suppose $\Delta \cup_f J$ is a generalized dunce hat. Since $J$ is a single circle, $f : \partial\Delta \to J$ is a taut map, and $f^{-1}(v)$ contains at least 3 points, then $\Delta \cup_f J$ is an unexposed taut one-relator presentation 2-complex which doesn't fall under the exceptional case described in the statement of the Theorem. Hence, the Theorem implies $\Delta \cup_f J$ is finitely unsplittable. ∎

Before embarking on the proof of the Theorem, we present evidence that weakening the hypotheses or strengthening the conclusions of the Theorem and the Corollary in various ways will vitiate them. Also we show that the conclusion of the Theorem can't be strengthened to rule out expressing an unexposed taut one-relator presentation 2-complex as the union of *three* proper subpolyhedra with trivial first homology.

First, we show that the scope of the Theorem can't be extended to include taut presentation 2-complexes with *naked* edges.

**Example 1.** Suppose $J = J_1 \vee J_2$ and $f : \partial\Delta \to J_1$ is a degree 3 covering map. (Hence, $f^{-1}(J_1 - \{v\})$ has at three components.) Then $Q = (\Delta \cup_f J_1) \vee J_2$ is a taut presentation 2-complexes with no free edges, but $J_2$ is a naked edge. Let $K$ and $B$ be subarcs of $J_2$ such that $v \in \text{int}(K)$, $v \notin B$ and $K \cup B = J_2$. Let $A = (\Delta \cup_f J_1) \cup K$. Then $A$ and $B$ are proper subpolyhedra of $Q$, $A \cup B = Q$, $H_1(A) \approx H_1(\Delta \cup_f J_1) \approx \mathbb{Z}_3$ and $H_1(B) = 0$. Hence, $Q$ is finitely splittable.

Second, we show that the Theorem becomes false if taut presentation 2-complexes are allowed to have *free* edges.

**Proposition 1.** If $Q = \Delta \cup_f J$ is a taut one-relator presentation 2-complex in which $J_1$ is a free edge, then $Q$ can be expressed as the union of two collapsible proper subpolyhedra. Hence, $Q$ is finitely splittable.

**Proof.** Let $J' = J_2 \vee J_3 \vee \ldots \vee J_n$. Then $A = f^{-1}(J') = \partial\Delta - f^{-1}(\text{int}(J_1))$ is an arc in $\partial\Delta$. Regard $\Delta$ as the product $A \times [0,1]$ where $A$ is identified with $A \times \{0\}$. Then the projection $A \times [0,1] \to A \times \{0\}$ can be regarded as a retraction $r : \Delta \to A$ with the property that for each subpolyhedron B of A, $r^{-1}(B)$ collapses to B. Conjugate r by $q^{-1}$ to obtain a retraction $s : Q \to J'$ with the property that for each subpolyhedron B of $J'$, $s^{-1}(B)$ collapses to B. For each i between 2 and n, express $J_i$ as the union of two subarcs $K_i$ and $L_i$ such that $v \in \partial K_i = \partial L_i = K_i \cap L_i$. Let $K = \cup_{i=2}^{n} K_j$ and $L = \cup_{i=2}^{n} L_j$. K and L are each cones with vertex v and $K \cup L = J'$. Hence, K and L are collapsible. Since $s^{-1}(K)$ and $s^{-1}(L)$ collapse to K and L, respectively, then $s^{-1}(K)$ and $s^{-1}(L)$ are each collapsible. Since $L \not\subset s^{-1}(K)$ and $K \not\subset s^{-1}(L)$, then $s^{-1}(K)$ and $s^{-1}(L)$ are proper subpolyhedra of Q. Since $s^{-1}(K) \cup s^{-1}(L) = Q$, then we have expressed Q as the union of two collapsible proper subpolyhedra. Hence, Q is finitely splittable. ∎

Third, we examine the three cases that arise when J is a single simple closed curve, $f : \partial \Delta \to J$ is a taut map and $f^{-1}(\{v\})$ has 1 or 2 elements.

**Example 3.** If $f^{-1}(\{v\})$ has 1 element, then $f : \partial \Delta \to J$ is a homeomorphism and $\Delta \cup_f J$ is a disk which can be expressed as the union of two proper subdisks. (Alternatively, observe that $\Delta \cup_f J$ has a free edge and apply Proposition 1.)

**Example 4.** Suppose $f^{-1}(\{v\})$ has 2 elements and $f : \partial \Delta \to J$ is a null-homotopic taut map. This is the exceptional case of the Theorem. Let $Q = \Delta \cup_f J$. Then the closures of the components of $\partial \Delta - f^{-1}(\{v\})$ are two arcs $A_1$ and $A_2$ in $\partial \Delta$ such that $A_1 \cup A_2 = \partial \Delta$ and $A_1 \cap A_2 = \partial A_1 = \partial A_2 = f^{-1}(\{v\})$. Furthermore, there is a homeomorphism $h : A_1 \to A_2$ such that $f \mid A_1 = f \circ h$. Then the quotient map $q : \Delta \to Q$ factors as $q = q_2 \circ q_1$ where $q_1 : \Delta \to \Delta/h$ and $q_2 : \Delta/h \to Q$ are quotient maps, $\Delta/h$ is the quotient space obtained from $\Delta$ by identifying x with $h(x)$ for each $x \in A_1$, and $q_2$ identifies the two distinct points of the set $q_1(f^{-1}(\{v\}))$. Note that $\Delta/h$ is a 2-sphere. Thus, Q is a 2-sphere with a pair of points identified which we denote $\mathcal{S}^2/\{2\ points\}$. Let $x \in int(A_1)$ and let B be an arc in $\Delta$ joining x to $h(x)$ such that $int(B) \subset int(\Delta)$. Then $\Delta$ is the union of two proper subdisks $D_1$ and $D_2$ such that $D_1 \cap D_2 = B$. Observe that $q(D_1) = q_2(q_1(D_1))$ and $q(D_2) = q_2(q_1(D_2))$ are proper subdisks of Q whose union is such Q. (Also $q(D_1) \cap q(D_2) = \partial q(D_1) \cup \{v\} = \partial q(D_2) \cup \{v\}$ where $\{v\} = int(q(D_1)) \cap int(q(D_2))$.) Thus, in this exceptional case, Q is finitely splittable.

**Example 5.** Suppose $f^{-1}(\{v\})$ has 2 elements and $f : \partial \Delta \to J$ is taut map which is not null-homotopic. Then $f : \partial \Delta \to J$ must be a degree two covering map, and $\Delta \cup_f J$ is a projective plane. Since the Theorem holds in this situation, then the projective plane $P^2$ must be finitely unsplittable. Here is a self-contained proof of this assertion. Suppose $P^2$ is the union of two proper subpolyhedra A and B with finite first homology groups. We make A connected without changing its first homology by joining distinct components of A by arcs. We do the same for B. Then we thicken A and B by passing to their regular neighborhoods and, thereby, can assume that A and B are compact connected surfaces with boundary in $P^2$. Then A and B collapse to connected 1-dimensional polyhedra K and L with the same first homology. Since the first homology of a 1-dimensional polyhedron is a direct sum of $\mathbb{Z}$'s, then $H_1(K) = 0 = H_1(L)$.

Hence, K and L must be trees. Thus, A and B are collapsible and must therefore be disks. But the only closed surface that can be expressed as the union of two disks is the 2-sphere, not $P^2$.

Fourth we explore whether the Theorem can be strengthened to conclude that an unexposed taut one-relator presentation 2-complex Q can't be represented as the union of two proper subpolyhedra A and B such that the inclusions of A and B into Q induce the zero map on first homology. In fact we show that this proposed strengthening is false for many generalized dunce hats including the original dunce hat.

**Example 6.** Consider a generalized dunce hat $Q = \Delta \cup_f J$ where J is a single simple closed curve and f : $\partial\Delta \to J$ is a taut map that induces an isomorphism between fundamental groups. (The original dunce hat has this property.) A Mayer-Vietoris sequence argument shows $H_1(Q) = 0$. Let D be a disk in int($\Delta$), let A = q(D) and B = q($\Delta$ – int(D)) where q : $\Delta \to Q$ is the quotient map. Then Q = A $\cup$ B, A is collapsible because it is a disk, and the inclusion of B into Q induces the zero map on first homology because $H_1(Q) = 0$.

If $\Delta \cup_f J$ is a generalized dunce hat, then f : $\partial\Delta \to J$ is a embedding on a neighborhood of every point of $\partial\Delta - f^{-1}(\{v\})$. Thus, f has only one *singular value* in J, namely v. We now explore the effect of forming a quotient space $\Delta \cup_f J$ for which the attaching map f : $\partial\Delta \to J$ is allowed to have two singular values in J. We will see that the Theorem fails for such quotient spaces. We will exhibit a situation in which f has two singular values in J and the resulting quotient space is splittable into two collapsible subpolyhedra that intersect in a collapsible subpolyhedron. This phenomenon plays a central role in the second author's 2014 Ph.D. thesis [S], allowing him to construct uncountably many topologically distinct contractible open 4-manifolds that are splittable in the sense of Gabai.

**Example 7.** Let $v_0$ and $v_1$ be distinct points of a simple closed curve J, and let $J_1$ and $J_2$ be the closures of the two components of J – { $v_0$, $v_1$ }. Let $\Delta$ be a disk and let $w_0, w_1, \ldots, w_5$ be six distinct points of $\partial\Delta$ listed in cyclic order. Label the closures of the components of $\partial\Delta - \{ w_0, w_1, \ldots, w_5 \}$ by $E_1, E_2, \ldots, E_6$ so that for $1 \le i \le 6$, $\partial E_i = \{ w_{i-1}, w_i \}$ (where $w_6 = w_0$). Let f : $\partial\Delta \to J$ be a map such that $f^{-1}(\{v_0\}) = \{ w_0, w_2, w_4 \}$, $f^{-1}(\{v_1\}) = \{ w_1, w_3, w_5 \}$, f maps each of $E_1, E_2$, and $E_3$ homeomorphically onto $J_1$, and f maps each of $E_4, E_5$, and $E_6$ homeomorphically onto $J_2$. (This description of f fixes f up to isotopies that preserve f$|\{ w_0, w_1, \ldots, w_5 \}$.) The quotient space $Q = \Delta \cup_f J$ is called the *Jester's hat* and originates in the second author's 2014 Ph.D. thesis where it is shown that the Jester's hat is the spine of a compact contractible 4-manifold whose interior is splittable in the sense of Gabai. The fact that the 4-manifold is splittable follows directly from the fact that the Jester's hat can be expressed as the union of two collapsible proper subpolyhedra A and B that intersect in a collapsible subpolyhedron. We sketch the construction of A and B. Let q : $\Delta \to Q$ denote the quotient map. Let L be an arc in $\Delta$ such that $\partial L = \{ w_0, w_3 \}$ and int(L) $\subset$ int($\Delta$). Then $\Delta$ is the union of two disks $D_1$ and $D_2$ where $\partial D_1 = E_1 \cup E_2 \cup E_3 \cup L$, $\partial D_2 = E_4 \cup E_5 \cup E_6 \cup L$ and $D_1 \cap D_2 = L$. Let $A = q(D_1)$ and $B = q(D_2)$. Hence, $Q = A \cup B$. A and B are easily seen to be collapsible. Indeed, the quotient map q transports a collapse of $D_1$ onto $E_1 \cup E_2 \cup E_3$ to a collapse of A onto the arc $J_1$, and $J_1$ collapses to a point. Similarly, B collapses to the arc $J_2$ which collapses to a point. Also $A \cap B = q(L)$ is an arc and is, thereby, collapsible.

We end this section by proving a proposition which tells us that the conclusion of the Theorem can't be strengthened to rule out the possibility of expressing an unexposed taut one-relator presentation 2-complex as the union of three proper subpolyhedra with trivial first homology.

**Proposition 2.** Every unexposed taut one-relator presentation 2-complex is the union of three proper collapsible subpolyhedra.

**Proof.** Let $Q = \Delta \cup_f J$ be an unexposed taut one-relator presentation 2-complex, where $J = J_1 \vee J_2 \vee \ldots \vee J_n$ is a wedge of n circles with wedge point v. Let $q : \Delta \to Q$ denote the quotient map. For each i between 1 and n, express $J_i$ as the union of two subarcs $K_i$ and $L_i$ such that $v \in \partial K_i = \partial L_i = K_i \cap L_i$. Let $K = \bigcup_{i=1}^{n} K_i$ and $L = \bigcup_{i=2}^{n} L_i$. K and L are each cones with vertex v and $K \cup L = J$. Let A and B be regular neighborhoods of K and L, respectively, in Q. Then A and B collapse to K and L, respectively, which in turn collapse to v. Thus A and B are collapsible. Since $q(\partial \Delta) = J = K \cup L \subset int(A) \cup int(B)$, then $q^{-1}(int(A)) \cup q^{-1}(int(B))$ is a neighborhood of $\partial \Delta$ in $\Delta$. Hence, there is a piecewise linear disk D in $int(\Delta)$ such that $\Delta - (q^{-1}(int(A)) \cup q^{-1}(int(B))) \subset D$. Let $C = q(D)$. Then C is a disk in Q. It follows that A, B and C are collapsible subpolyhedra of Q that cover Q. ∎

**3. Auxiliary facts about components of submanifolds and finite graphs.**
Here we record two propositions that will be used in the proof of the Theorem.

**Notation.** We need notations that distinguish between the two types of *interiors* that can be associated with a manifold that is the subset of a larger topological space. If A is a subset of a topological space X, let *clx(A), intx(A)* and *frx(A)* denote the *closure, interior* and *frontier* of A in X, respectively. Omit the subscript *X* in contexts where no ambiguity can arise from this omission. If M is a topological manifold, then from this point on, let $in(M)$ and $\partial M$ denote the manifold interior and the manifold boundary of M, respectively.

**Definition.** If $M \subset Q$ are manifolds and $dim(M) = n$, we call M a *nice submanifold* of Q if $in(M) \subset in(Q)$ and if $\partial M \cap \partial Q$ and $\partial M - in(\partial M \cap \partial Q)$ are (n − 1)-manifolds.

Suppose T is a combinatorial triangulation of an n-manifold M (i.e., the link of every vertex of T is piecewise linearly homeomorphic to either a piecewise linear (n – 1)-sphere or a piecewise linear (n – 1)-ball.)  Let K be a subcomplex of T.  Then regular neighborhood theory ([Br], Theorem 3.6, page 231) implies that the simplicial neighborhood of K in a second derived subdivision T″ of T is an n-dimensional nice submanifold of M.

**Notation.**  If X is a topological space, let *Comp(X)* denote the set of all components of X.

**Proposition 3.**  Suppose M ⊂ Q are manifolds such that $in$(M) ⊂ $in$(Q).  Then the functions C ↦ cl$_M$(C) : *Comp*(M ∩ $in$(Q)) → *Comp*(M) and D ↦ D ∩ $in$(Q) : *Comp*(M) → *Comp*(M ∩ $in$(Q)) are well defined bijections which are each other's inverses.

**Proof.**  It suffices to prove that these two functions are well defined and that each is a left inverse of the other.

We will repeatedly use the following facts.

*i)*   If A is a connected subset of a topological space X, and A ⊂ B ⊂ cl$_X$(A), then B is connected and cl$_X$(B) = cl$_X$(A).

*ii)*   Every component of a topological space X is a closed subset of X.

*iii)*   If C is a component of an m-manifold M, then C is an open subset of M (because M is locally connected) and, thus, C is an m-manifold.

*iv)*   If M is a non-empty connected m-manifold, then so is $in$(M) (because any arc in M joining two points of $in$(M) can be pushed into $in$(M)).

*v)*   If a manifold D is a closed subset of a topological space X, then cl$_X$($in$(D)) = D (because $in$(D) is a dense subset of D).

*vi)*   If M ⊂ Q are manifolds of the same dimension, then *Invariance of Domain* implies that $in$(M) ⊂ $in$(Q).  Hence, in this situation, M ∩ ∂Q = ∂M ∩ ∂Q.

First we prove that the function C ↦ cl$_M$(C) : *Comp*(M ∩ $in$(Q)) → *Comp*(M) is well defined, and that the function D ↦ D ∩ $in$(Q) is a left inverse of C ↦ cl$_M$(C).  Let C ∈ *Comp*(M ∩ $in$(Q)).  Since C is a connected subset of M, then there is a D ∈ *Comp*(M) such that C ⊂ D.  We will prove cl$_M$(C) = D and D ∩ $in$(Q) = C.  Since D is a closed subset of M (by fact *ii)*), then cl$_M$(C) ⊂ D.  Let m = dim(M).  Since M ∩ $in$(Q) is an open subset of M, then M ∩ $in$(Q) is an m-manifold.  Hence, C is an m-manifold (by fact *iii)*).  Also D is an m-manifold (by fact *iii)*).  Since C ⊂ D ⊂ M, then $in$(C) ⊂ $in$(D) ⊂ $in$(M) (by fact *vi*)). Since *in*(M) ⊂ *in*(Q), then *in*(D) ⊂ M ∩ *in*(Q). Also *in*(D) is connected (by fact *iv*)) and intersects *in*(C). Since C is the component of M ∩ *in*(Q) that contains *in*(C), then *in*(D) ⊂ C. Therefore (by fact *v*)), D = cl$_M$(*in*(D)) ⊂ cl$_M$(C). Hence, cl$_M$(C) = D. This proves C ↦ cl$_M$(C) : *Comp*(M ∩ *in*(Q)) → *Comp*(M) is a well defined function. Clearly, C ⊂ D ∩ *in*(Q). Since *in*(D) is connected and *in*(D) ⊂ D ∩ *in*(Q) ⊂ D = cl$_M$(*in*(D)) (by fact *v*)), then D ∩ *in*(Q) is connected (by fact *i*)). Since *in*(C) ⊂ D ∩ *in*(Q) ⊂ M ∩ *in*(Q) and C is the component of M ∩ *in*(Q) that contains *in*(C), then D ∩ *in*(Q) ⊂ C. We conclude that D ∩ *in*(Q) = C. This proves that the function D ↦ D ∩ *in*(Q) is a left inverse of C ↦ cl$_M$(C).

To complete the proof of Proposition 3, we must show that the function D ↦ D ∩ *in*(Q) : *Comp*(M) → *Comp*(M ∩ *in*(Q)) is well defined, and that the function C ↦ cl$_M$(C) is a left inverse of D ↦ D ∩ *in*(Q). Let D ∈ *Comp*(M). Then D and M are both m-manifolds (by fact *iii*)). Hence, *in*(D) ⊂ *in*(M) (by fact *vi*)). Since D ⊂ M and *in*(M) ⊂ *in*(Q), then *in*(D) ⊂ D ∩ *in*(Q) ⊂ D = cl$_M$(*in*(D)) (by fact *v*)). Since *in*(D) is connected (by fact *iv*)), then D ∩ *in*(Q) is connected (by fact *i*)). D is a closed and open subset of M (by facts *ii*) and *iii*)). Hence, D ∩ *in*(Q) is a closed and open subset of M ∩ *in*(Q). It follows that D ∩ *in*(Q) is a component of M ∩ *in*(Q). This proves D ↦ D ∩ *in*(Q) : *Comp*(M) → *Comp*(M ∩ *in*(Q)) is a well defined function. Let C = D ∩ *in*(Q). Then C ∈ *Comp*(M ∩ *in*(Q)). Since D is the component of M such that C ⊂ D, then the argument in the previous paragraph shows cl$_M$(C) = D. This proves C ↦ cl$_M$(C) is a left inverse of D ↦ D ∩ *in*(Q). ∎

**Definition.** By a *finite graph* we mean a 1-dimensional CW complex with finitely many 0-cells and 1-cells. The *vertices* of the finite graph are its 0-cells and the *edges* are its 1-cells. form a set with the discrete topology. An *edge* of the graph is obtained by attaching an *arc* (a homeomorph of [0,1]) by an *attaching map* from the boundary of the arc (a 2-point set) to the set of vertices. Thus, two different types of edges as possible. If the attaching map is injective, then the resulting edge is an arc joining two distinct vertices. Alternatively, if the attaching map is a constant map, then the resulting edge is a simple closed curve containing a single vertex; such an edge is called a *loop.* This definition allows for the existence of two distinct edges which are arcs joining the same pair of vertices; such edges are sometimes called *multiple edges.* If Γ is a finite graph, then the *polyhedron underlying Γ*, denoted |Γ|, is the union of the vertices and edges of Γ. We give |Γ| the *Whitehead topology:* a subset U of |Γ| is open if and only if for every edge E of Γ, U ∩ E is a relatively open subset of E (which is topologically an arc or a simple closed curve). Because Γ is a finite graph, the Whitehead topology on |Γ| is metrizable. (Assign each edge of Γ length 1 and declare the distance between two points of |Γ| to be the length of the shortest path joining them in |Γ|.)

**Definition.**  A topological space T is a *finite tree* if T is the polyhedron underlying a finite graph and T contains no simple closed curves.

**Definition.**  The *degree* of a vertex v of a finite graph Γ, denoted *deg(v,Γ)* or simply *deg(v)* (if the identity of Γ is clear from context) is the number of edges of Γ that emanate from v where an edge of Γ which is a loop with vertex v contributes 2 to the degree of v.  A vertex v of a finite graph Γ is *isolated* if deg(v,Γ) = 0 and is an *endpoint* of Γ if deg(v,Γ) = 1.  Let $\mathcal{V}(\Gamma)$ and $\mathcal{E}(\Gamma)$ denote the set of vertices and the set of edges, respectively, of the finite graph Γ.

**Definition.**  Suppose Γ is a finite graph.  Let $\Gamma_{\neq 2}$ be the maximal subcomplex of Γ that has no vertices of degree 2.  Then $|\Gamma| - |\Gamma_{\neq 2}|$ is a 1-manifold without boundary.  Let $\mathcal{C}$ denote the set of all closures of components of $|\Gamma| - |\Gamma_{\neq 2}|$.  Then for each C ∈ $\mathcal{C}$, one of the following three cases occurs.

*i)*  C is an arc and for v ∈ $\mathcal{V}(\Gamma)$ ∩ C, deg(v,Γ) = 2 if and only if v ∈ int(C).

*ii)*  C is a simple closed curve and there is exactly one element v of $\mathcal{V}(\Gamma)$ ∩ C such that deg(v,Γ) ≠ 2.

*iii)*  C is a simple closed curve and a component of $|\Gamma|$ and, thus, deg(v,Γ) = 2 for all v ∈ $\mathcal{V}(\Gamma)$ ∩ C.

Suppose that no component of $|\Gamma|$ is a simple closed curve.  Then only cases *i)* and *ii)* can occur.  In this situation, we form a finite graph Γ* from Γ by *consolidating* each element of $\mathcal{C}$ into a single edge of Γ*, and we call Γ* the *consolidation* of Γ.  Observe that there is a subcomplex $\Gamma_{\neq 2}$* of Γ* such that $\Gamma_{\neq 2}$* contains all the vertices of Γ*, and there is a homeomorphism $h_{\Gamma}$* : $|\Gamma|$ → $|\Gamma$*$|$ such that $h_{\Gamma}$*($|\Gamma_{\neq 2}|$) = $|\Gamma_{\neq 2}$*$|$, $h_{\Gamma}$* $\big|$ $|\Gamma_{\neq 2}|$ : $|\Gamma_{\neq 2}|$ → $|\Gamma_{\neq 2}$*$|$ is induced by a simplicial isomorphism from $\Gamma_{\neq 2}$ to $\Gamma_{\neq 2}$*, and for each C ∈ $\mathcal{C}$, $h_{\Gamma}$*(C) is an edge of Γ* that is not an edge of $\Gamma_{\neq 2}$*.  It follows that for each v ∈ $\mathcal{V}(\Gamma_{\neq 2})$, deg(v,Γ) = deg($h_{\Gamma}$*(v),Γ*), and Γ* has no vertices of degree 2.  Thus, for each k ≥ 0, if k ≠ 2, then Γ and Γ* have the same number of vertices of degree k.  Note that the process of creating Γ* by consolidation from Γ may create new loops and multiple edges.  However, this process can't produce simple closed curves in $|\Gamma$*$|$ unless they already exist in $|\Gamma|$.  Thus, if each component of $|\Gamma|$ is a finite tree, then so is each component of $|\Gamma$*$|$.  Note that if each component of $|\Gamma|$ is a finite tree, then only case *i)* can occur, and distinct elements of $\mathcal{C}$ are arcs that are either disjoint or have at most one endpoint in common.

**Lemma 1.**  Suppose that Σ is a finite graph such that each component of $|\Sigma|$ is a finite tree.  If each vertex of Σ that is not an endpoint has degree ≥ 3, then the number of endpoints of Σ is strictly greater than the number of non-endpoint vertices of Σ.

**Proof.**  It suffices to consider the case in which Σ is a finite tree.  Let V and E denote the number of vertices and edges of Σ, respectively.  Since Σ is a finite tree, its Euler characteristic is 1. Thus, $V - E = 1$. Let $V_\partial$ and $V_\iota$ denote the number of endpoints and non-endpoint vertices of T, respectively. Then $V = V_\partial + V_\iota$. Since $2E$ counts the vertices of $\Sigma$ with multiplicity, then $2E \geq V_\partial + 3V_\iota$. Hence,

$$2V_\partial + 2V_\iota \; = \; 2V \; = \; 2E + 2 \; \geq \; V_\partial + 3V_\iota + 2.$$

Thus, $V_\partial \geq V_\iota + 2$. ∎

**Lemma 2.** If $\Gamma$ is a finite graph in which each vertex has degree $\geq 3$, then the number of vertices of $\Gamma$ that lie in simple closed curves in $|\Gamma|$ is strictly greater than the number of vertices that do not.

**Proof.** Let $\Sigma$ denote the subcomplex of $\Gamma$ consisting of all the edges of $\Gamma$ that do not lie in simple closed curves in $|\Gamma|$ together with the vertices of these edges. Since every vertex of $\Sigma$ lies in an edge of $\Sigma$, then $\Sigma$ has no isolated vertices. Also, clearly, $|\Sigma|$ contains no simple closed curves. Thus, each component of $|\Sigma|$ is a finite tree.

Let $\Sigma^*$ denote the consolidation of $\Sigma$ as defined above. Then the remarks in the paragraph containing the definition of *consolidation* imply that each component of $|\Sigma^*|$ is a finite tree, $\Sigma^*$ has no vertices of degree 2, and for $k \geq 0$, if $k \neq 2$, then $\Sigma$ and $\Sigma^*$ have the same number of vertices of degree k. Thus, $\Sigma^*$ has no isolated vertices. Hence, Lemma 1 implies that $\Sigma^*$ has more endpoints than non-endpoint vertices.

Let $V_s$ denote the number of vertices of $\Gamma$ that lie in simple closed curves in $|\Gamma|$, and let $V_o$ denote the number of other vertices of $\Gamma$. We must prove $V_o < V_s$. Let $V_\partial$ denote the number of endpoints of $\Sigma$, and let $V_{\geq 3}$ denote the number of vertices of $\Sigma$ of degree $\geq 3$. Finally, let $V_\partial^*$ and $V_\iota^*$ denote the number of endpoints and non-endpoint vertices of $\Sigma^*$, respectively. Then the previous lemma implies $V_\iota^* < V_\partial^*$. The remarks of the previous paragraph imply $V_{\geq 3} = V_\iota^*$ and $V_\partial^* = V_\partial$.

We assert that $V_0 \leq V_{\geq 3}$. Indeed, suppose that v is a vertex of $\Gamma$ that lies in no simple closed curve in $|\Gamma|$. Then every edge of $\Gamma$ that contains v must belong to $\Sigma$. This implies that v lies in $\Sigma$ and $\deg(v,\Sigma) = \deg(v,\Gamma) \geq 3$. Hence, $V_0 \leq V_{\geq 3}$.

We also assert that $V_\partial \leq V_s$. Suppose v is an endpoint of $\Sigma$. Then $\deg(v,\Sigma) = 1 < 3 \leq \deg(v,\Gamma)$. Hence, some edge of $\Gamma$ that contains v doesn't belong to $\Sigma$. It follows that this edge must lie in a simple closed curve in $|\Gamma|$. So v must lie in a simple closed curve in $|\Gamma|$. This proves $V_\partial \leq V_s$.

Combining the preceding information yields: $V_0 \leq V_{\geq 3} = V_\iota^* < V_\partial^* = V_\partial \leq V_s$. ∎

**Proposition 4.** If $\Gamma$ is a finite graph in which each vertex has degree $\geq 2$ and $\alpha : \mathcal{V}(\Gamma) \to \mathcal{V}(\Gamma)$ is an involution of $\mathcal{V}(\Gamma)$ without fixed points that preserves the degree of each vertex, then there is a vertex v of $\Gamma$ such that both v and $\alpha(v)$ lie in simple closed curves in $|\Gamma|$. Furthermore, if $\Gamma$ has a vertex of degree $> 2$, then we can choose v to be a vertex of $\Gamma$ of degree $> 2$.

**Proof.** First observe that if every vertex of $\Gamma$ has degree 2, then $|\Gamma|$ is a closed 1-manifold. Hence, every vertex of $\Gamma$ lies on a simple closed curve in $|\Gamma|$ – finishing the proof. So we might as well assume that $\Gamma$ has vertices of degree > 2. Let $\beta$ denote the restriction of $\alpha$ to $\{ v \in \mathcal{V}(\Gamma) : \deg(v) > 2 \}$. It suffices to find a vertex $v$ of $\Gamma$ of degree > 2 such that both $v$ and $\beta(v)$ lie in simple closed curves in $|\Gamma|$. Next, we might as well delete from $\Gamma$ all components of $|\Gamma|$ that contain no vertices of degree > 2; in other words, assume that each component of $|\Gamma|$ contains a vertex of degree > 2.

Form the consolidation $\Gamma^*$ of $\Gamma$ as defined above. Since each vertex of $\Gamma$ has degree $\geq 2$, then each vertex of $\Gamma^*$ has degree $\geq 3$. Recall that there is a homeomorphism $h_\Gamma{}^* : |\Gamma| \to |\Gamma^*|$ such that $(h_\Gamma{}^*)^{-1}$ sends each vertex of $\Gamma^*$ to a vertex of $\Gamma$ of the same degree ($\geq 3$). Since $\beta$ preserves degrees of vertices, it induces an involution $\beta^* = h_\Gamma{}^* \circ \beta \circ (h_\Gamma{}^*)^{-1} | \mathcal{V}(\Gamma^*) : \mathcal{V}(\Gamma^*) \to \mathcal{V}(\Gamma^*)$ without fixed points that also preserves degrees of vertices. Thus, for some positive integer $n$, we can write $\mathcal{V}(\Gamma^*)$ as the disjoint union of two n-element sets $\mathcal{V}_1$ and $\mathcal{V}_2$ such that $\beta^*(\mathcal{V}_1) = \mathcal{V}_2$. Let $\mathcal{S}^*$ denote the set of vertices of $\Gamma^*$ that lie in simple closed curves in $|\Gamma^*|$.

We will prove that there is a $v^* \in \mathcal{S}^* \cap \mathcal{V}_1$ such that $\beta^*(v^*) \in \mathcal{S}^* \cap \mathcal{V}_2$. If not, then $\beta^*(\mathcal{S}^* \cap \mathcal{V}_1) \cap (\mathcal{S}^* \cap \mathcal{V}_2) = \emptyset$. For each finite set A, let #A denote the number of elements in A. Since $\beta^*(\mathcal{S}^* \cap \mathcal{V}_1) \subset \mathcal{V}_2$, then $\#\beta^*(\mathcal{S}^* \cap \mathcal{V}_1) + \#(\mathcal{S}^* \cap \mathcal{V}_2) \leq \#\mathcal{V}_2 = n$. Since $\beta^*$ is injective, then $\#(\mathcal{S}^* \cap \mathcal{V}_1) + \#(\mathcal{S}^* \cap \mathcal{V}_2) \leq n$. Since $\mathcal{S}^*$ is the disjoint union of $\mathcal{S}^* \cap \mathcal{V}_1$ and $\mathcal{S}^* \cap \mathcal{V}_2$, then $\#\mathcal{S}^* \leq n$. But the previous lemma implies that $\#\mathcal{S}^* > n$. We conclude that there is a $v^* \in \mathcal{S}^* \cap \mathcal{V}_1$ such that $\beta^*(v^*) \in \mathcal{S}^* \cap \mathcal{V}_2$. Thus, $v^*$ is a vertex of $\Gamma^*$ of degree $\geq 3$ such that both $v^*$ and $\beta^*(v^*)$ lie in simple closed curves in $|\Gamma^*|$. Let $v = (h_\Gamma{}^*)^{-1}(v^*)$. Then $v$ is a vertex of $\Gamma$ of degree $\geq 3$, and the homeomorphism $(h_\Gamma{}^*)^{-1} : |\Gamma^*| \to |\Gamma|$ sends the simple closed curves in $|\Gamma^*|$ that contain $v^*$ and $\beta^*(v^*)$ to simple curves in $|\Gamma|$ that contain $v$ and $(h_\Gamma{}^*)^{-1}(\beta^*(v^*)) = \beta(v)$. ∎

**4. Proof of the Theorem.** Let $Q = \Delta \cup_f J$ be an unexposed taut one-relator presentation 2-complex where $\Delta$ is a 2-dimensional disk, $J = J_1 \vee J_2 \vee \ldots \vee J_n$ is a wedge of n circles with vertex v, and $f : \partial\Delta \to J$ is a taut map. Because Q has no exposed edges, then $f^{-1}(J_i - \{v\})$ must have at least 2 components for each i between 1 and n.

Next we discuss the *link of v in Q* because it plays an important role in the upcoming proof.

**Definition.** Suppose P is a polyhedron, $v \in P$, and K is a triangulation of P such that v is a vertex of K. Then the *star of v in K*, denoted $St(v,K)$, is the subcomplex of K consisting of all simplices of K that contain v and their faces, and the *link of v in K*, denoted $\mathcal{L}k(v,K)$, is the subcomplex of $St(v,K)$ consisting of all simplices of $St(v,K)$ that don't contain v. Let $St(v,P)$ and $Lk(v,P)$ denote the polyhedra underlying $St(v,K)$ and $\mathcal{L}k(v,K)$, respectively. ($St(v,P)$ and $Lk(v,P)$ are the unions of all the simplices of K that belong to $St(v,K)$ and $Lk(v,K)$, respectively.)   Then St(v,P) and Lk(v,P) are piecewise linear invariants of the pair (P,v).  ([Br], Theorem 2.1 and its proof, page 225.) Moreover, St(v,P) and Lk(v,P) naturally embed as subpolyhedra of P.

Suppose $Q = \Delta \cup_f J$ is an unexposed taut one-relator presentation 2-complex. We will now describe specific embeddings of St(v,Q) and Lk(v,Q) in Q.  For each i between 1 and n, let $E_i$ be a small arc in $J_i$ containing v in its interior.  We identify St(v,$J_i$) and Lk(v,$J_i$) with $E_i$ and $\partial E_i$, respectively.  Observe that $\cup_{i=1}^{n} E_i$ is a cone in J with vertex v and base $\cup_{i=1}^{n} \partial E_i$.  We identify St(v,J) with $\cup_{i=1}^{n} E_i$ and Lk(v,J) with $\cup_{i=1}^{n} \partial E_i$.  Observe that $f^{-1}(\cup_{i=1}^{n} E_i)$ is the union of a finite pairwise disjoint collection $\mathcal{L}$ of arcs in $\partial\Delta$ such that $f(\cup\mathcal{L}) = \cup_{i=1}^{n} E_i$.  (Obviously, $f(\cup\mathcal{L}) \subset \cup_{i=1}^{n} E_i$.  Equality is a consequence of the fact that Q has no naked edges.)  For each $L \in \mathcal{L}$, there is a point $u_L \in in(L)$ such that $\{u_L\} = f^{-1}(\{v\}) \cap L$.  For each $L \in \mathcal{L}$, one of two possibilities occurs.  Either $f|L : L \to \cup_{i=1}^{n} E_i$ is an embedding which maps the two components of $L - \{u_L\}$ homeomorphically onto two distinct components of $(\cup_{i=1}^{n} E_i) - \{v\}$, or $f|L$ "folds" L at the interior point $u_L$ and maps each of the two components of $L - \{u_L\}$ homeomorphically onto the same component of $(\cup_{i=1}^{n} E_i) - \{v\}$.  Let $\{ D_L : L \in \mathcal{L} \}$ be a finite pairwise disjoint collection of disks in $\Delta$ such that $D_L \cap (\partial\Delta) = L$ for each $L \in \mathcal{L}$.  For each $L \in \mathcal{L}$, let $L' = \partial D_L - in(L)$.  Then $\{ L' : L \in \mathcal{L} \}$ is a finite pairwise disjoint collection of arcs in $\Delta$ such that for each $L \in \mathcal{L}$, $\partial L' = \partial L$ and $in(L') \subset in(\Delta)$.  We identify St(v,Q) and Lk(v,Q) with the subsets $q(\cup_{L \in \mathcal{L}} D_L)$ and $q(\cup_{L \in \mathcal{L}} L')$ of Q, respectively.  Then St(v,Q) $\cap$ J = St(v,J) = $\cup_{i=1}^{n} E_i$ and Lk(v,Q) $\cap$ J = Lk(v,J) = $\cup_{i=1}^{n} \partial E_i$, and in particular for $1 \le i \le n$,  St(v,Q) $\cap J_i = E_i$ and Lk(v,Q) $\cap J_i = \partial E_i$. Observe that Lk(v,Q) is a 1-dimensional polyhedron.  We regard Lk(v,Q) as the polyhedron underlying a finite graph $\Lambda(v,Q)$ whose set of vertices is Lk(v,J) = $\cup_{i=1}^{n} \partial E_i = \cup_{L \in \mathcal{L}} q(\partial L')$  and whose set of edges is $\{ q(L') : L \in \mathcal{L} \}$.  Notice that for each $L \in \mathcal{L}$, $q(\partial L')$ is a set consisting of either two vertices or one vertex of $\Lambda(v,Q)$; in the former case, $q(L')$ is an arc joining two vertices, and in the latter case, $q(L_j')$ is a simple closed curve containing a single vertex.

Let w be a vertex of $\Lambda = \Lambda(v,Q)$.  Then $w \in \partial E_i \subset J_i - \{v\}$ for some i between 1 and n.  $f^{-1}(\{w\})$ is a finite set with one element in each component of $f^{-1}(J_i - \{v\})$.  For each $L \in \mathcal{L}$, if an endpoint of L lies $f^{-1}(\{w\})$, then $q(L')$ is an edge of $\Lambda$ emanating from w. Moreover, $q(L')$ contributes 1 to deg(w,$\Lambda$) if only one endpoint of L lies in $f^{-1}(\{w\})$, while $q(L')$ contributes 2 to deg(w,$\Lambda$) if both endpoints of L lie in $f^{-1}(\{w\})$.  (In the former situation, $q(L')$ is an arc with one endpoint at w, while in the latter situation, $q(L')$ is a simple closed curve containing w.)  Moreover, each point of $f^{-1}(\{w\})$ is an endpoint of L' for exactly one element L of $\mathcal{L}$.  Thus, deg(w,$\Lambda$) equals the number of elements of $f^{-1}(\{w\})$ which, in turn, equals the number of components of $f^{-1}(J_i - \{v\})$.  Since Q has no exposed edges, then $f^{-1}(J_i - \{v\})$ has at least 2 components, and, thus, deg(w,$\Lambda$) $\ge 2$. It follows that if $f^{-1}(J_i - \{v\})$ has exactly 2 components for $1 \le i \le n$, then deg(w,$\Lambda$) = 2 for each vertex w of $\Lambda$.  In that case, each component of Lk(v,Q) is a simple closed curve. We note that if $\partial E_i = \{ w, w' \}$, then deg(w,$\Lambda$) and deg(w',$\Lambda$) both equal to the number of components of $f^{-1}(J_i - \{v\})$ and, hence, are equal to each other.  Therefore, a degree preserving involution $\alpha : \mathcal{V}(\Lambda) \to \mathcal{V}(\Lambda)$ without fixed points is defined by having $\alpha$ interchange the points of $\partial E_i$ for $1 \le i \le n$. Proposition 4 now implies that, in the case that some vertex of $Lk(v,Q)$ has degree $\ge 3$, then there is an i between 1 and n such that both points of $\partial E_i$ lie in simple closed curves in $Lk(v,Q)$ and have degree $\ge 3$. Thus, if $f^{-1}(J_i - \{v\})$ has at least 3 components for some i between 1 and n, then there is an i between 1 and n such that both points of $\partial E_i = Lk(v,Q) \cap J_i$ lie in simple closed curves in $Lk(v,Q)$, and $f^{-1}(J_i - \{v\})$ has at least 3 components.

Assume Q is finitely splittable. Then Q is the union of two proper subpolyhedra A and B each of which has a finite first homology group.

At this point, the proof bifurcates into two cases.

**Case 1.  $f^{-1}(J_i - \{v\})$ has exactly two components for $1 \le i \le n$.**

**Case 2.  $f^{-1}(J_i - \{v\})$ has at least 3 components for some i between 1 and n.**

**The proof of Case 1.**  Assume $v \in A$. Although we didn't initially suppose that A and B are connected, we now make them connected. We join distinct components of A by arcs (one fewer arcs than the number of components). According to the reduced Mayer-Vietoris sequence, this process doesn't alter A's first homology. We do the same for B.

In this case, each point of $J_i - \{v\}$ has a neighborhood in Q homeomorphic to $\mathbb{R}^2$. Thus, $Q - \{v\}$ is an open 2-manifold. Also $Lk(v,K)$ is the union of a finite pairwise disjoint collection $\{ S_1, S_2, \ldots , S_k \}$ of simple closed curves, and $St(v,K) = (v \star S_1) \cup (v \star S_2) \cup \ldots \cup (v \star S_k)$, where each $v \star S_i$ is a cone with vertex v and base $S_i$ that is embedded in Q so that $(v \star S_i) \cap (v \star S_j) = \{v\}$ for $1 \le i < j \le k$.

The quotient map $q : \Delta \to Q$ factors as $q = q_2 \circ q_1$ where $q_1 : \Delta \to Q'$, $q_2 : Q' \to Q$ and $Q'$ is the closed 2-manifold obtained from Q by replacing each cone $v \star S_i$ by a cone $w_i \star S_i$ for $1 \le i \le k$ so that $(w_i \star S_i) \cap (w_j \star S_j) = \emptyset$ for $1 \le i < j \le k$. $q_2 : Q' \to Q$ is the identity on each $S_i$ and outside $\bigcup_{i=1}^{k} (w_i \star S_i)$, $q_2(w_i) = v$ and $q_2$ maps each arc in the cone structure of $w_i \star S_i$ to the corresponding arc in the cone structure of $v \star S_i$ for $1 \le i \le k$. $q_1 : \Delta \to Q'$ is the quotient map that identifies the two components of the $f^{-1}(J_i - \{v\})$'s but identifies two elements of $f^{-1}(\{v\})$ only if it is forced by the identifications of the two components of the $f^{-1}(J_i - \{v\})$'s. (There is a one-to-one correspondence between the points of $q_1(f^{-1}(\{v\}))$ and the components of $Lk(v,K)$.)

Let T be a triangulation of Q so that subcomplexes of T triangulate $\{v\}$, each $J_i$, A, B, and so that there is a triangulation $q_2^{-1}(T)$ of $Q'$ that makes $q_2 : Q' \to Q$ a simplicial map from $q_2^{-1}(T)$ to T. Let $T''$ be a second derived subdivision of T. Then there is a second derived subdivision $q_2^{-1}(T)''$ of $q_2^{-1}(T)$ so that $q_2 : Q' \to Q$ is simplicial from $q_2^{-1}(T)''$ to $T''$. We enlarge A and B by replacing them by their simplicial neighborhood with respect to T″.  This process replaces $q_2^{-1}(A)$ and $q_2^{-1}(B)$ by their simplicial neighborhood with respect to $q^{-1}(T)$″.  Then according to the remark following the definition of *nice submanifold* in section 3, $q_2^{-1}(A)$ and $q_2^{-1}(B)$ are compact 2-dimensional submanifolds of Q′.  Thus, $q_2^{-1}(A)$ and $q_2^{-1}(B)$ have non-empty boundaries.

Since $v \in A$ and A is "thickened" to a simplicial neighborhood of itself, then $q_2^{-1}(\{v\}) \subset in(q_2^{-1}(A))$.

We assert that every component $q_2^{-1}(A)$ must contain a point of $q_2^{-1}(\{v\}) = \{ w_1, w_2, \ldots , w_k \}$.  For if C is a component of $q_2^{-1}(A)$ that is disjoint form $q_2^{-1}(\{v\})$, then A is the union of the two disjoint closed subsets $q_2(C)$ and $q_2(q_2^{-1}(A) - C)$.  This contradicts the fact that A is connected.

No component of $q_2^{-1}(A)$ can contain two distinct elements of $q_2^{-1}(\{v\}) = \{ w_1, w_2, \ldots , w_k \}$.  For if C is a component of $q_2^{-1}(A)$ that contains $w_1$ and $w_2$, then C contains an arc E with endpoints $w_1$ and $w_2$ that misses $w_3, w_4, \ldots , w_k$.  There is a retraction $r´ : q_2^{-1}(A) \to E$ such that $r´(q_2^{-1}(A) - C) = \{w_1\}$.  If we conjugate $r´$ by $(q_2 \,|\, q_2^{-1}(A))^{-1}$, we obtain a retraction r of A onto the simple closed curve $q_2(E)$.  Hence, $id_{q_2(E)} = r \circ i$ where $i : q_2(E) \to A$ denotes inclusion.  This factorization induces the factorization of id : $H_1(q_2(E)) \to H_1(q_2(E))$ as the composition of $i_\star : H_1(q_2(E)) \to H_1(A)$ and $r_\star : H_1(A) \to H_1(q_2(E))$.  Since $H_1(q_2(E)) \approx \mathbb{Z}$ and $H_1(A)$ is a finite group, this is impossible.  It follows that for each component C of $q_2^{-1}(A)$, $q_2 \,|\, C : C \to Q$ is an embedding.

Each component of $q_2^{-1}(A)$ is a disk.  For suppose C is a component of $q_2^{-1}(A)$ such that $C \cap q_2^{-1}(\{v\}) = \{w_1\}$ and let $D = q_2^{-1}(A) - C$.  Then $A = q_2(C) \cup q_2(D)$ and $q_2(C) \cap q_2(D) = \{v\}$.  Hence, $H_1(C,\{w_1\}) \approx H_1(q_2(C),\{v\}) \approx H_1(A,q_2(D)$ (by excision).  $q_2(D)$ is connected because every component of $q_2(D)$ contains the point v.  The exactness of the sequence $H_1(A) \to H_1(A,q_2(D)) \to \tilde{H}_0(q_2(D))$ implies $H_1(A,q_2(D))$ is a finite group.  Hence, $H_1(C,\{w_1\})$ is a finite group.  The exactness of the sequence, $H_1(\{w_1\}) \to H_1(C) \to H_1(C,\{w_1\}) \to \tilde{H}_0(\{w_1\})$ implies $H_1(C)$ is a finite group.  Since C is a compact connected 2-manifold with boundary, it collapses to a compact connected 1-dimensional polyhedron K with the same first homology group.  Since the first homology group of a compact 1-dimensional polyhedron is a direct sum of $\mathbb{Z}$'s, then $H_1(K) = 0$.  Thus, K is a tree.  Therefore, C is collapsible and must be a disk.

We can apply the same arguments to B and $q_2^{-1}(B)$ to show that each component of $q_2^{-1}(B)$ is a disk.  (If $v \in B$, then the argument of the previous paragraph applies.  If $v \notin B$, then B is a compact connected 2-manifold with non-empty boundary and finite first homology group.  Then, as argued in the previous paragraph, each component of B and of $q_2^{-1}(B)$ is a disk.)  Also, as we showed for the components of $q_2^{-1}(A)$, no component of $q_2^{-1}(B)$ can contain two points of $\{ w_1, w_2, \ldots , w_k \}$.  We have $Q´ = q_2^{-1}(A) \cup q_2^{-1}(B)$.  We can enlarge $q_2^{-1}(A)$ to a single disk by connecting distinct components of $q_2^{-1}(A)$ by bands. Similarly we can enlarge $q_2^{-1}(B)$ to a single disk. Thus, $Q'$ is the union of 2 disks. Since the only closed 2-manifold that can be expressed as the union of 2 disks is a 2-sphere, then $Q'$ must be a 2-sphere.

The map $q_1 : \Delta \rightarrow Q'$ imposes a cell-structure on $Q'$. This cell-structure has one 2-cell $q_1(\Delta)$, n 1-cells $q_1(cl(f^{-1}(J_i - \{v\})))$ for $1 \leq i \leq n$, and k vertices $q_1(f^{-1}(\{v\})) = \{ w_1, w_2, \ldots , w_k \}$. Therefore, the Euler characteristic of $Q'$ is $k - n + 1$. Consequently, $k - n + 1 = 2$. So $k = n + 1$. Thus, $k \geq 2$.

We will prove $k \leq 2$. Assume $k \geq 3$. Let $D_i$ be the component of $q_2^{-1}(A)$ that contains $w_i$ for $1 \leq i \leq k$. Then $q_2^{-1}(A) = \bigcup_{i=1}^{k} D_i$. Therefore, $Q' - (\bigcup_{i=1}^{k} D_i) \subset q_2^{-1}(B)$. Since $Q' - (\bigcup_{i=1}^{k} D_i)$ is connected, then it is contained in a single disk component C of $q_2^{-1}(B)$. C can contain at most one $w_i$. Therefore, after reindexing, we can suppose $w_1$ and $w_2 \notin C$. Push $\partial D_1$ slightly out of $D_1$ to a simple closed curve L in C that separates $w_1$ and $w_2$ in the 2-sphere $Q'$. There is a retraction of $Q' - \{ w_1, w_2 \}$ onto L that restricts to a retraction r of C onto L. Then $id_L = r \circ j$ where $j : L \rightarrow C$ denotes inclusion. This composition induces a composition This composition induces the factorization of id : $H_1(L) \rightarrow H_1(L)$ as the composition of $j_* : H_1(L) \rightarrow H_1(C)$ and $r_* : H_1(C) \rightarrow H_1(L)$. Since $H_1(q_2(E)) \approx \mathbb{Z}$ and $H_1(C) = \{0\}$, this is impossible. We have proved that $k \leq 2$. Hence, $k = 2$ and $n = k - 1 = 1$. Thus, J is a single simple closed curve. Since $f^{-1}(J - \{v\})$ has exactly 2 components, we are in one of the two situations described in Examples 4 and 5 in Section 2. In Example 4, f : $\partial \Delta \rightarrow J$ is null-homotopic and Q is $\mathbb{S}^2 / \{2 \text{ points}\}$, the single exceptional case mentioned in the statement of the Theorem. In Example 5, f : $\partial \Delta \rightarrow J$ is a degree 2 covering map and Q is a projective plane which is finitely unsplittable and, therefore, satisfies the conclusion of the Theorem.

This completes the proof of Case 1. ∎

**The proof of Case 2.** Since $f^{-1}(J_i - \{v\})$ has at least 3 components for some i between 1 and n, then after reindexing, we can assume that $f^{-1}(J_1 - \{v\})$ has at least 3 components and both points of $Lk(v,Q) \cap J_1$ lie in simple closed curves in $Lk(v,Q)$.

To guide the reader through the upcoming argument, we briefly describe the five steps that comprise it.

*i)*   Find a point $b \in J_1 - B$ such that $b \neq v$.

*ii)*   Make A connected by adding arcs and then thicken A to a regular neighborhood of itself. This makes each $A \cap J_i$ a 1-dimensional submanifold of $J_i$, $A - J$ a 2-dimensional submanifold of $P - J$ and $q^{-1}(A)$ a nice 2-dimensional submanifold of $\Delta$.

*iii)*   Let $\mathcal{E}$ denote the set of components of $A \cap J$ and let $\mathcal{F}$ denote the set of components of $A - J$. Construct a graph $\Gamma$ in A with a vertex $w_E \in E$ for each E $\in \mathcal{E}$, a vertex $w_F \in F$ for each F $\in \mathcal{F}$, and an edge joining $w_E$ to $w_F$ in cl(F) if and only if $E \cap cl(F) \neq \varnothing$. Then prove that A strong deformation retracts onto l$\Gamma$l, l$\Gamma$l is a tree, and cl(F) is a disk for each $F \in \mathcal{F}$ and A is the union of { cl(F) : $F \in \mathcal{F}$ }. Thus, A is a *tree of disks.*

*iv)* Prove **Lemma 4:** If a disk G in Q intersects J in an arc $H \subset J_1 - \{v\}$ such that $\partial H \subset \partial G$ and $b \in in(H) \subset in(G)$, then the inclusion map $\partial G \to Q - \{b\}$ induces a monomorphism on first homology.

*v)* Let $E_b$ be the unique element of $\mathcal{E}$ that contains b. Construct a subcomplex $\Gamma^*$ of $\Gamma$ that contains the vertex $w_{E_b}$ and is maximal in $\Gamma$ with respect to the property that for each vertex $w_E$ of $\Gamma^*$, except possibly $w_{E_b}$, deg($w_E, \Gamma^*$) = 2. The nature of $\Gamma^*$ at $w_{E_b}$ breaks into three cases depending on the local topology of A near b. The upshot is that the union A* of { cl(F) : $F \in \mathcal{F}$ and $w_F$ is a vertex of $\Gamma^*$ } is a **disk** containing b in its interior with the property that $\partial A^* \subset fr(A) \subset B$. It follows that the inclusion induced map $H_1(\partial A^*) \to H_1(Q - \{b\})$ factors through the finite group $H_1(B)$. The latter statement contradicts the conclusion of Lemma 4 because we can choose a disk G in $in(A^*)$ satisfying the hypotheses of Lemma 4.

If you understand these five steps completely, read no further. For everyone else, we provide a detailed proof.

We begin by proving B can't contain $J_1$.

**Lemma 3.** If $Q = \Delta \cup_f J$ is an unexposed taut one-relator presentation 2-complex and B is a proper subpolyhedron of Q with finite first homology group, then B can't contain $J_i$ for any i between 1 and n.

**Proof.** Suppose B contains $J_1$. Since no $J_i$ is a naked edge of Q for $1 \leq i \leq n$, then Q − J is a dense subset of Q. Hence, every non-empty open subset of Q must contain points of Q − J. Since Q − B is a non-empty open subset of Q, then there is a point $x \in (Q − B) − J$. There is a point $y \in in(\Delta)$ such that $q^{-1}(\{x\}) = \{y\}$. Consider a retraction of $\Delta − \{y\}$ onto $\partial\Delta$. We conjugate this retraction by $q^{-1}$ to obtain a retraction $r_1$ of Q − {x} onto J. Let $r_2$ be the retraction of J onto $J_1$ which, for $2 \leq i \leq n$, maps $J_i$ to v. Then $r_2 \circ r_1 | B$ is a retraction r of B onto $J_1$. Then $r \circ i = id_{J_1}$ where $i : J_1 \to B$ denotes the inclusion. Hence, the identity map $H_1(J_1) \to H_1(J_1)$ factors through the finite group $H_1(B)$. Since $H_1(J_1) \approx \mathbb{Z}$, this is impossible. Thus, B can't contain $J_1$. ∎

Observe that the argument used in the proof of Lemma 3 also shows neither A nor B can contain $J_i$ for any i between 1 and n.

Since $J_1 − B$ is a non-empty open subset of $J_1$, we can choose a point $b \in J_1 − B$ such that $b \neq v$. Since $A \cup B = Q$, then necessarily $b \in A$.

Although we did not initially suppose that A is connected, we now make A connected as we did in the proof of Case 1.

Let T be a triangulation of Q so that subcomplexes of T triangulate {v}, each $J_i$, A, B and {b}, and so that there is a triangulation $q^{-1}(T)$ of $\Delta$ that makes $q : \Delta \rightarrow Q$ a simplicial map from $q^{-1}(T)$ to T. Let $T''$ be a second derived subdivision of T. Then there is a second derived subdivision $q^{-1}(T)''$ of $q^{-1}(T)$ so that $q : \Delta \rightarrow D$ is simplicial from $q^{-1}(T)''$ to $T''$. We enlarge A by replacing it by its simplicial neighborhood with respect to $T''$. This process replaces $q^{-1}(A)$ by its simplicial neighborhood with respect to $q^{-1}(T)''$. Then according to the remark following the definition of *nice submanifold* in section 3, $A \cap J_i$ is now a 1-dimensional submanifold of $J_i$ for $1 \leq i \leq n$, and $q^{-1}(A)$ is now a nice 2-dimensional submanifold of $\Delta$. Thus, $q^{-1}(A) \cap \partial\Delta$ is a 1-dimensional submanifold of $\partial\Delta$, and $q^{-1}(A) \cap in(\Delta)$ is a 2-dimensional submanifold of $in(\Delta)$. Since q maps $in(\Delta)$ homeomorphically onto Q − J, then it follows that A − J is a 2-dimensional submanifold of Q − J. Also $b \in in(A \cap J_1)$. Furthermore, if $v \in A$, then $v \in in(A \cap J_i)$ for $1 \leq i \leq n$.

We devote the next three paragraphs to describing the topology of A at points of $A \cap J$. Let $\mathcal{E}$ denote the set of components of $A \cap J$, and for each i between 1 and n, let $\mathcal{E}_i$ denote the set of components of $A \cap J_i$. Since $A \cap J_i$ is a 1-dimensional submanifold of $J_i$, then each element of $\mathcal{E}_i$ is an arc in $J_i$. Let $E \in \mathcal{E}_i$. Since v is a vertex of the triangulation T of Q and E arises by taking a simplicial neighborhood of $A \cap J_i$ in a second derived subdivision $T''$ of T, then $v \in E$ implies $v \in in(E)$. Similarly, if i = 1 and b $\in E$, then $b \in in(E)$. Now consider an element E of $\mathcal{E}$. If $v \notin E$, then E is an element of some $\mathcal{E}_i$ and, hence E is an arc in $J_i − \{v\}$. (Then, as previously noted, if i = 1 and $b \in E$, then $b \in in(E)$.) On the other hand, if $v \in E$, then $E = \cup_{i=1}^{n} E_i$ where $E_i \in \mathcal{E}_i$ and $v \in in(E_i)$ for $1 \leq i \leq n$. Thus, E is a cone in J with vertex v and base $\cup_{i=1}^{n} \partial E_i$. Therefore, E is a regular neighborhood of v in J. Again either $b \notin E$ or $b \in in(E_1)$. Since $f = q \,|\, \partial\Delta : \partial\Delta \rightarrow$ J maps each component of $\partial\Delta − q^{-1}(\{v\})$ homeomorphically onto $J_i − \{v\}$ for some i, $1 \leq i \leq$ n, it follows that for each $E \in \mathcal{E}$, $q^{-1}(E)$ is the union of a finite pairwise disjoint collection of arcs in $\partial\Delta$ such that if $v \in E$, then $q^{-1}(\{v\}) \subset in(q^{-1}(E))$ and if $b \in E$, then $q^{-1}(\{b\}) \subset in(q^{-1}(E))$.

Consider an $E \in \mathcal{E}$. Then $q^{-1}(E)$ is the union of a finite pairwise disjoint collection $\mathcal{D}$ of arcs in $\partial\Delta$. First focus on the case in which $v \notin E$. In this situation q maps each D $\in \mathcal{D}$ homeomorphically onto E. Let N be a regular neighborhood of E in A. Then $q^{-1}(N)$ is a regular neighborhood of $q^{-1}(E)$ in $q^{-1}(A)$. Since $q^{-1}(A)$ is a nice 2-dimensional submanifold of $\Delta$, it follows that $q^{-1}(N)$ is the union of a finite pairwise disjoint collection of disks { $N_D : D \in \mathcal{D}$ } such that $N_D \cap \partial\Delta = D$ for each $D \in \mathcal{D}$. Furthermore, q embeds each $N_D$ in N such that $q(N_{D_1}) \cap q(N_{D_2}) = E$ for distinct $D_1$ and $D_2 \in \mathcal{D}$. Thus, N is an *open book* with *pages* $q(N_D)$, $D \in \mathcal{D}$, and binding E.

Next consider the case in which $v \in E$. Recall that in this situation, $E = \bigcup_{i=1}^{n} E_i$ where $E_i \in \mathcal{E}_i$ and $v \in in(E_i)$ for $1 \le i \le n$. As we previously observed, $E$ is a cone in $J$ with vertex $v$ and base $\bigcup_{i=1}^{n} \partial E_i$. For each $i$ between 1 and n, express $E_i$ as the union of two subarcs $E_i{}'$ and $E_i{}''$ such that $E_i{}' \cap E_i{}'' = \{v\}$. Again $q^{-1}(E)$ is the union of a finite pairwise disjoint collection $\mathcal{D}$ of arcs in $\partial\Delta$. In this situation, for each $D \in \mathcal{D}$, there is a point $u_D \in in(D)$ such that $q^{-1}(\{v\}) \cap D = \{u_D\}$ and $D$ can be expressed as the union of two subarcs $D'$ and $D''$ such that $D' \cap D'' = \{u_D\}$. Then for each $D \in \mathcal{D}$, $q$ maps each of $D'$ and $D''$ homeomorphically onto one of the arcs $E_1{}'$, $E_1{}''$, $E_2{}'$, $E_2{}''$, ... , $E_n{}'$, $E_n{}''$. Furthermore, if $q(D') \ne q(D'')$, then $q$ embeds $D$ in $E$. However, if $q(D') = q(D'')$, then $q \mid D$ has a single critical point at $u_D$. Let $\mathcal{D}_1 = \{ D \in \mathcal{D} : q(D') \ne q(D'') \}$ and let $\mathcal{D}_2 = \{ D \in \mathcal{D} : q(D') = q(D'') \}$. Then $\mathcal{D}$ is the disjoint union of $\mathcal{D}_1$ and $\mathcal{D}_2$. Again consider a regular neighborhood $N$ of $E$ in $A$. As before, $q^{-1}(N)$ is a regular neighborhood of $q^{-1}(E)$ in $q^{-1}(A)$. Since $q^{-1}(A)$ is a nice 2-dimensional submanifold of $\Delta$, it follows that $q^{-1}(N)$ is the union of a finite pairwise disjoint collection of disks $\{ N_D : D \in \mathcal{D} \}$ such that $N_D \cap \partial\Delta = D$ for each $D \in \mathcal{D}$. Furthermore, $q$ embeds each $N_D - D$ in $N - E$ such that $q(N_{D_1} - D_1) \cap q(N_{D_2} - D_2) = \emptyset$ for distinct $D_1$ and $D_2 \in \mathcal{D}$. Thus, $N$ is an *open pop-up book* with pages $q(N_D)$, $D \in \mathcal{D}$, and binding $E$. Moreover, for $D \in \mathcal{D}_1$, $q$ embeds $N_D$ in $N$ and the page $q(N_D)$ is a disk attached to $E$ along the arc $q(D') \cup q(D'')$ in its boundary; whereas for $D \in \mathcal{D}_2$, $q \mid N_D$ is not an embedding and the page $q(N_D)$ is a "cone" with vertex $v$ attached to $E$ along the arc $q(D') = q(D'') \in \{ E_1{}', E_1{}'', E_2{}', E_2{}'', ... , E_n{}', E_n{}'' \}$ which joins the vertex $v$ of this cone to a point in its base.

Next we argue that each component of $q^{-1}(A)$ is a disk. We first observe that every component of $q^{-1}(A)$ must intersect $\partial\Delta$. For suppose $C$ is a component of $q^{-1}(A)$ that is disjoint from $\partial\Delta$. Then $C$ and $q^{-1}(A) - C$ are disjoint closed subsets of $\Delta$, and $q(C)$ and $q(q^{-1}(A) - C)$ are disjoint closed subsets of $Q$ whose union is $A$. This contradicts our earlier supposition that $A$ is connected.

Now suppose $C$ is a component of $q^{-1}(A)$ that is not a disk. Then $C$ is a disk with holes and there is a component $K$ of $\partial C$ that lies $in(\Delta)$ and bounds a disk $F$ in $in(\Delta)$. Since every component of $q^{-1}(A)$ intersects $\partial\Delta$, then $in(F)$ can't contain a component of $q^{-1}(A)$. Thus, $in(F) \cap q^{-1}(A) = \emptyset$. Since $\Delta - in(F)$ is an annulus with boundary components $\partial\Delta$ and $K$, and since $q^{-1}(A) \cap \partial\Delta$ is a proper subset of $\partial\Delta$, then there is a retraction of $\Delta - in(F)$ onto $K$ that maps $q^{-1}(A) \cap \partial\Delta$ to a single point $z \in K$. We restrict this retraction to $C$ to obtain a retraction of $C$ onto $K$ that maps $C \cap \partial\Delta$ to $z$. We conjugate this retraction by $(q \mid C)^{-1}$ to obtain a retraction $r_0$ of $q(C)$ onto $q(K)$ that maps $q(C) \cap J$ to $q(z)$. Let $D = q^{-1}(A) - C$. Then $C$ and $D$ are disjoint closed subsets of $\Delta$ such that $q(C) \cup q(D) = A$ and $q(C) \cap q(D) \subset J$. Since $r_0$ maps $q(C) \cap J$ to $q(z)$, then $r_0$ extends to a retraction $r$ of $A$ onto $q(K)$ which sends $q(D)$ to $q(z)$. Hence, $r \circ j = \mathrm{id}_{q(K)}$ where $j : q(K) \to A$ denotes the inclusion. Therefore, the identity map $H_1(q(K)) \to$ $H_1(q(K))$ factors through the finite group $H_1(A)$. Since $H_1(q(K)) \approx \mathbb{Z}$, this is impossible. Thus, each component of $q^{-1}(A)$ must be a disk.

The next step of the proof is to construct a graph $\Gamma$ in A which is a spine of A in the sense that there is a strong deformation retraction of A onto $\Gamma$. Recall that $\mathcal{E}$ denotes the set of components of $A \cap J$, and each element of $\mathcal{E}$ is either an arc in $J - \{v\}$ or a regular neighborhood of v in J (which happens to be a cone with vertex v). Let $\mathcal{F}$ denote the set of components of $A - J$. Then each $F \in \mathcal{F}$ is a 2-dimensional submanifold of $Q - J$, and cl(F) is a compact subpolyhedron of Q such that, for each component D of cl(F) $\cap$ J, there is an element E of $\mathcal{E}$ such that $D \subset E$. Furthermore, if v $\notin$ E, then D = E. However, if v $\in$ E (and E is a cone in J with vertex v and base $\beta(E)$), then either D is an arc in E joining two distinct points of $\beta(E)$, or D is an arc in E joining v to a point of $\beta(E)$.

$\Gamma$ has two types of vertices which we choose according to the following rules.

*i*)  For each $E \in \mathcal{E}$: if v $\notin$ E, choose a point $w_E \in in(E)$, and if v $\in$ E, choose $w_E = v$.

*ii*)  For each $F \in \mathcal{F}$, choose a point $w_F \in in(F)$.

Then $\{ w_E : E \in \mathcal{E} \} \cup \{ w_F : F \in \mathcal{F} \}$ is the set of vertices of $\Gamma$. Since q maps $q^{-1}(A) \cap in(\Delta)$ homeomorphically onto $A - J$, then $\{ q^{-1}(F) : F \in \mathcal{F} \}$ is the set of components of $q^{-1}(A) \cap in(\Delta)$. For each $F \in \mathcal{F}$, let $P(F) = cl(q^{-1}(F))$. Since $q^{-1}(A)$ is a nice submanifold of $\Delta$, then Proposition 3 implies that $\{ P(F) : F \in \mathcal{F} \}$ is the collection of components of $q^{-1}(A)$, and $P(F) \cap in(\Delta) = q^{-1}(F)$ for each $F \in \mathcal{F}$. Hence, each P(F) is a nice 2-dimensional submanifold of $\Delta$ and, in particular, as we argued in a previous paragraph, P(F) is a disk. Furthermore, $q^{-1}(w_F) \in in(P(F))$. Let $B(F) = P(F) \cap q^{-1}(\{ w_E : E \in \mathcal{E} \})$. Then B(F) is a finite subset of $P(F) \cap \partial\Delta \subset \partial P(F)$. Let K(F) be a cone in P(F) with vertex $q^{-1}(w_F)$ and base B(F) such that $K(F) \cap \partial P(F) = B(F)$. Then q(K(F)) is a 1-dimensional polyhedron in cl(F). Let $\Gamma_F$ be the graph whose underlying polyhedron is q(K(F)), whose vertex set is $\{ w_F \} \cup q(B(F))$ and whose edges are the arcs in q(K(F)) joining $w_F$ to a point of q(B(F)). Let $\Gamma = \cup_{F \in \mathcal{F}} \Gamma_F$. Then $\Gamma$ is a graph in A. (In fact $\Gamma$ is a *bipartite* graph whose only edges join a vertex of the set $\{ w_E : E \in \mathcal{E} \}$ to a vertex of the set $\{ w_F : F \in \mathcal{F} \}$.) Also observe that because $q^{-1}(A)$ is a nice compact 2-dimensional submanifold of $\Delta$, then each component of $q^{-1}(A) \cap \partial\Delta$ must be an arc which lies in the boundary of a unique component of $q^{-1}(A)$. Thus, each component of $q^{-1}(A) \cap \partial\Delta$ lies in $\partial P(F)$ for a unique $F \in \mathcal{F}$.

There is a strong deformation retraction of A onto $|\Gamma|$, the polyhedron underlying $\Gamma$, which is constructed in several steps. Begin with a strong deformation retraction of A $\cap$ J onto $\{ w_E : E \in \mathcal{E} \}$. (Since each component E of A $\cap$ J is either an arc with $w_E \in in(E)$ or a cone with vertex $w_E = v$, then this strong deformation retraction clearly exists.) Each component D of    $q^{-1}(A \cap J)$ is an arc which is mapped by q into a component E of A ∩ J in a fashion that allows the strong deformation retraction of E onto $w_E$ to be lifted to a strong deformation retraction of D onto the one-point set D ∩ $q^{-1}(\{w_E\})$ via conjugation by q│D. Hence, the strong deformation retraction of A ∩ J onto { $w_E$ : E ∈ $\mathcal{E}$ } can be lifted (via conjugation by q) to a strong deformation retraction of $q^{-1}(A ∩ J)$ onto $q^{-1}(\{ w_E : E ∈ \mathcal{E} \})$.    For each F ∈ $\mathcal{F}$, this strong deformation retraction restricts to a strong deformation retraction of P(F) ∩ ∂Δ onto B(F) which, in turn, extends to a strong deformation retraction of P(F) onto K(F).  This strong deformation retraction pushes down (via conjugation by $q^{-1}$) to a strong deformation retraction of cl(F) onto q(K(F)) which extends the original strong deformation retraction of cl(F) ∩ J onto { $w_E$ : E ∈ $\mathcal{E}$ } ∩ cl(F).  The union of these strong deformation retractions is a strong deformation retraction of A onto IΓI.

The end of the strong deformation retraction of A onto IΓI is a retraction map r : A → IΓI.  Then r∘σ = $id_{IΓI}$ where σ : IΓI → A denotes the inclusion.  Hence, the identity map $H_1(IΓI) → H_1(IΓI)$ factors through the finite group $H_1(A)$.  Thus, $H_1(IΓI)$ must be a finite group.  Since Γ is a graph, then $H_1(IΓI)$ is a direct sum of $\mathbb{Z}$'s.  Hence, $H_1(IΓI) = 0$ and IΓI must be a tree.

Our next goal is to show that cl(F) is a disk for each F ∈ $\mathcal{F}$.  Let F ∈ $\mathcal{F}$.  q maps the disk P(F) onto cl(F), and q maps P(F) − ∂Δ homeomorphically onto cl(F) − J.  It remains to analyze the behavior of q on the components of P(F) ∩ ∂Δ.  q can't identify two distinct points of B(F), because if it did, this would create a simple closed curve in IΓI, contradicting the fact that IΓI is a tree.  So if x ∈ $q^{-1}(\{w_E\})$ and x′ ∈ $q^{-1}(\{w_{E'}\})$ are distinct points of B(F) where E and E′ ∈ $\mathcal{E}$, then $w_E ≠ w_{E'}$ and, hence, E ≠ E′.  Therefore, q must map distinct components of P(F) ∩ ∂Δ into distinct components of A ∩ J.  (Said another way: if E ∈ $\mathcal{E}$ and cl(F) ∩ E ≠ ∅, then P(F) contains exactly one of the components of $q^{-1}(E)$ and is disjoint from all the others.)  Suppose E ∈ $\mathcal{E}$ and cl(F) ∩ E ≠ ∅.  Then there is a unique component E* of P(F) ∩ ∂Δ which q maps into E.  There are three possible cases.

*i)*  v ∉ E.  In this case, q maps E* homeomorphically onto E.

*ii)*  v ∈ E and q embeds E* in E.  In this case, E is a cone in J with vertex v and base β(E) and q maps E* homeomorphically onto an arc in E joining two distinct points of β(E).

*iii)*  v ∈ E and q│E* : E* → E has a critical point.  In this case, there is a unique point v* ∈ $in(E*)$ such that q(v*) = v, and q maps each component of E* − {v*} homeomorphically onto a single component of E − {v}.  (In other words, q maps E* onto an arc in E joining v to a point of β(E) *folding* E at v*.)  The effect that q has on P(F) is to *zip up* P(F) along the arc E* ⊂ ∂P(F), taking $in(E*)$ to a half open arc that lies in $in(cl(F))$.

q maps P(F) onto cl(F) and q maps P(F) − ∂Δ homeomorphically onto cl(F) − J.  In the event that case *iii)* doesn't occur, then q maps each component of P(F) ∩ ∂Δ homeomorphically onto a corresponding component of cl(F) ∩ J. In this situation q maps P(F) homeomorphically onto cl(F), making cl(F) a disk. Also note that in this situation, the components of cl(F) ∩ J = q(P(F) ∩ Δ) are arcs in ∂cl(F). In the event that case *iii)* occurs (for a unique E ∈ $\mathcal{E}$ for which v ∈ E and cl(F) ∩ E ≠ ∅), then q maps all but one of the components of P(F) ∩ ∂Δ homeomorphically onto a corresponding component of cl(F) ∩ J, and q *zips up* the remaining component of P(F) ∩ ∂Δ as described above. In this situation, although q│P(F) : P(F) → cl(F) is not a homeomorphism, nonetheless q│P(F) does not alter the topological type of P(F). Hence, cl(F) is also a disk in the event that case *iii)* occurs.

The preceding description reveals:

for E ∈ $\mathcal{E}$ and F ∈ $\mathcal{F}$, $w_E$ and $w_F$ are joined by an edge in Γ if and only if E ∩ cl(F) ≠ ∅.

Indeed, an edge of Γ joining $w_E$ to $w_F$ is the image under q of an arc in P(F). Thus, $w_E$ ∈ E ∩ q(P(F)) = E ∩ cl(F). On the other hand, if E ∩ cl(F) ≠ ∅, then the way that $w_E$ is chosen guarantees that $w_E$ ∈ cl(F) = q(P(F)). Thus, P(F) ∩ $q^{-1}$({$w_E$}) contains a point x, the cone K(F) contains an arc joining its vertex $q^{-1}(w_F)$ to x, and q maps this arc to an edge of Γ that joins $w_F$ to $w_E$.

Since A is the union of the collection of disks { cl(F) : F ∈ $\mathcal{F}$ } and the intersection pattern of { cl(F) : F ∈ $\mathcal{F}$ } is encoded by Γ where lΓl is a tree, we call A a *tree of disks*.

Recall that the point b lies in $J_1$ − B and b ≠ v. Hence, b ∈ A and there is an $E_b$ ∈ $\mathcal{E}$ such that b ∈ $E_b$. Furthermore, if v ∉ $E_b$, then $E_b$ ∈ $\mathcal{E}_1$, $E_b$ is an arc and b ∈ *in*($E_b$). On the other hand, if v ∈ $E_b$, then $E_b$ is a cone in J with vertex v. Furthermore, in this situation, for each i between 1 and n, there is an arc $E_i$ ∈ $\mathcal{E}_i$ such that v ∈ *in*($E_i$), $E_i$ = $E_b$ ∩ $J_i$, $E_b$ = $\cup_{i=1}^{n} E_i$ and b ∈ *in*($E_1$). When v ∈ $E_b$, we identify $E_b$ with St(v,J) = St(v,Q) ∩ J. Then, our choice at the beginning of the proof of Case 2 guarantees that each of the two points of ∂$E_1$ = Lk(v,Q) ∩ $J_1$ lie in simple closed curves in Lk(v,Q). (The same choice insures that $f^{-1}(J_1 − \{v\})$ has at least 3 components.)

To complete the proof of the Theorem, we will need:

**Lemma 4.** If G is a disk in Q which intersects J in an arc H ⊂ $J_1$ − {v} such that ∂H ⊂ ∂G and b ∈ *in*(H) ⊂ *in*(G), then the inclusion map $j$ : ∂G → Q − {b} induces a monomorphism on first homology.

**Proof.** Let $G_1$ and $G_2$ be the closures of the two components of G − H. Then $G_1$ and $G_2$ are disks in Q such that $G_i$ ∩ J = H ⊂ ∂$G_i$ for i = 1, 2 and $G_1$ ∩ $G_2$ = H. For i = 1, 2, let $G_i$* = cl($q^{-1}(G_i − H)$). Then $G_i$* is a disk in Δ which q maps homeomorphically onto $G_i$. Also $H_i$* = $G_i$* ∩ ∂Δ is an arc in ∂Δ which is a component of $q^{-1}(H)$. *in*($H_1$*) and *in*($H_2$*) must be disjoint because any sufficiently small neighborhood N in Δ of a point of $in(H_1{}^*) \cap in(H_2{}^*)$ would have the property that N $\cap$ $in(\Delta)$ is contained in both $G_1{}^* - H_1{}^*$ and $G_2{}^* - H_2{}^*$. However, $G_1{}^* - H_1{}^*$ and $G_2{}^* - H_2{}^*$ are disjoint because q maps them to disjoint sets. Since $in(H_1{}^*) \cap in(H_2{}^*) = \emptyset$, then $H_1{}^* \neq H_2{}^*$. Since distinct components of $q^{-1}(H)$ must lie in distinct components of $q^{-1}(J_1 - \{v\})$, then $H_1{}^*$ and $H_2{}^*$ lie in distinct components of $q^{-1}(J_1 - \{v\})$. Since $q^{-1}(J_1 - \{v\})$ has at least 3 components, we can choose distinct components $M_1$, $M_2$ and $M_3$ of $q^{-1}(J_1 - \{v\})$ so that $H_1{}^* \subset M_1$ and $H_2{}^* \subset M_2$.

Form a cone T in $\Delta$ with base $\beta(T) = \partial\Delta - q^{-1}(J_1 - \{v\}) = q^{-1}(J_2 \vee J_3 \vee \ldots \vee J_n)$ and with vertex at a point $z \in in(\Delta)$ such that $T \cap \partial\Delta = \beta(T)$. There is a homeomorphism of $\Delta$ which fixes $\partial\Delta$ and pulls $G_1{}^*$ and $G_2{}^*$ off T. The inverse of this homeomorphism moves T off $G_1{}^* \cup G_2{}^*$. So we can assume T is disjoint from $G_1{}^* \cup G_2{}^*$. $H_1{}^*$ and $H_2{}^*$ lie in distinct components of $\Delta - T$. Hence, so do $G_1{}^*$ and $G_2{}^*$. For each $x \in \beta(T)$, let $T_x$ denote the arc in the cone structure of T joining z to x. For $1 \leq i \leq 3$, let $\{x_i, y_i\}$ be the endpoints of the open arc $M_i$ listed so that their cyclic order on $\partial\Delta$ is $x_1, y_1, x_2, y_2, x_3, y_3$ allowing for the possibility that $y_1 = x_2$, $y_2 = x_3$ and/or $y_3 = x_1$. Then $\partial\Delta - (M_1 \cup M_2 \cup M_3)$ has 3 components $M_4$, $M_5$ and $M_6$, each of which is either an arc or a one-point set listed so that $\{y_1, x_2\} \subset M_4$, $\{y_2, x_3\} \subset M_5$ and $\{y_3, x_1\} \subset M_6$. Thus, $M_1$, $M_2$, … , $M_6$ appear on $\partial\Delta$ in the cyclic order $M_1$, $M_4$, $M_2$, $M_5$, $M_3$, $M_6$.

Since each component of $q^{-1}(J_1 - \{v\})$ contains exactly one element of the set $q^{-1}(\{b\})$, then there is a unique retraction of $\partial\Delta - q^{-1}(\{b\})$ onto $\beta(T)$. This retraction extends to a retraction $r^*$ of $\Delta - q^{-1}(\{b\})$ onto T that maps $\partial G_i{}^* - in(H_i{}^*)$ homeomorphically onto $T_{x_i} \cup T_{y_i}$ for $i = 1, 2$. $r^*$ induces a retraction r of Q $- \{b\}$ onto q(T) (via conjugation by $q^{-1}$) that retracts J $- \{b\}$ to $J_2 \vee J_3 \vee \ldots \vee J_n$ and maps the open arc $\partial G_i - H$ homeomorphically onto the open arc $q(T_{x_i} \cup T_{y_i}) - \{v\}$ for $i = 1, 2$.

Let C denote the wedge of two circles $C_1$ and $C_2$ which intersect in the one-point set $\{c\}$. Let $s : q(T) \to C$ be a map which satisfies the following conditions: s maps $q(\{z\} \cup \beta(T))$ to $\{c\}$, s maps $q(in(T_x))$ homeomorphically onto $C_1 - \{c\}$ for $x \in \beta(T) \cap M_6$, s maps $q(in(T_x))$ to $\{c\}$ for $x \in \beta(T) \cap M_4$ and s maps $q(in(T_x))$ homeomorphically onto $C_2 - \{c\}$ for $x \in \beta(T) \cap M_5$. $H_1(C) \approx \mathbb{Z} \oplus \mathbb{Z}$ and there is a generating set $\{\gamma_1, \gamma_2\}$ for $H_1(C)$ such that the subgroup of $H_1(C)$ generated by $\gamma_i$ is the image of the inclusion induced map $H_1(C_i) \to H_1(C)$. Observe that for $i = 1,2$, $s \circ r \circ j \,|\, \partial G_i - in(H)$ is homotopic rel endpoints to a map which takes $\partial G_i - H$ homeomorphically onto $C_i - \{c\}$ and takes $\partial H$ to $\{c\}$. Thus, the induced map $s_* \circ r_* \circ j_* : H_1(\partial G) \to H_1(C)$ maps a generator of $H_1(\partial G) \approx \mathbb{Z}$ onto one of $\pm(\gamma_1 + \gamma_2)$ or $\pm(\gamma_1 - \gamma_2)$. Since these elements of $H_1(C)$ are of infinite order, then $s_* \circ r_* \circ j_*$ is a monomorphism. It follows that $j_* : H_1(\partial G) \to H_1(Q - \{b\})$ must be a monomorphism. ∎

Finally we come to the heart of the argument. Recall that $b \in J_1 - B \subset A \cap J_1$, $b \neq v$ and $E_b \in \mathcal{E}$ such that $b \in E_b$. We will break the remainder of the proof of Case 2 into three subcases. We will deal with the first subcase before we explain the second and third subcases.

**Subcase 2$i$): $v \notin E_b$.** In this subcase, $E_b \in \mathcal{E}_1$, $E_b$ is an arc and $b \in in(E_b)$.

Let $F_1, \ldots, F_p$ be the distinct elements of $\mathcal{F}$ such that $w_{F_1}, \ldots, w_{F_p}$ are all the vertices of $\Gamma$ that are joined to $w_{E_b}$ by an edge of $\Gamma$. We will prove $p \geq 3$. Since $E_b$ is a component of $A \cap J$, then each component of $q^{-1}(E_b)$ is a component of $q^{-1}(A) \cap \partial\Delta$. Since $\{ P(F) : F \in \mathcal{F} \}$ is the set of all components of $q^{-1}(A)$, and $P(F) \cap \partial\Delta \subset \partial P(F)$ for each $F \in \mathcal{F}$, then each component of $q^{-1}(E_b)$ is contained in $\partial P(F)$ for some $F \in \mathcal{F}$. When this occurs, $E_b \cap cl(F) \neq \emptyset$ and $w_{E_b}$ and $w_F$ are joined by an edge of $\Gamma$. Two different components of $q^{-1}(E_b)$ can't lie in $\partial P(F)$ for the same $F \in \mathcal{F}$, because then there would be two distinct edges of $\Gamma$ joining $w_{E_b}$ to $w_F$ contradicting the fact that $|\Gamma|$ is a tree. Thus, $p$ is at least as large as the number of distinct components of $q^{-1}(E_b)$. Since $v \notin E_b$, then each component of $f^{-1}(J_1 - \{v\})$ contains a distinct component of $q^{-1}(E_b)$. Our choice at the beginning of the proof of Case 2 implies that $f^{-1}(J_1 - \{v\})$ has at least 3 components. We conclude that $p \geq 3$.

Since $|\Gamma|$ is a tree and $w_{E_b}$ has $p$ edges emanating from it, then $|\Gamma| - \{ w_{E_b} \}$ has $p$ components. Let $\Gamma_1, \Gamma_2, \ldots, \Gamma_p$ be the $p$ subgraphs of $\Gamma$ such that $|\Gamma_1|, |\Gamma_2|, \ldots, |\Gamma_p|$ are the closures of the $p$ components of $|\Gamma| - \{ w_{E_b} \}$ labeled so that $w_{F_i}$ is a vertex of $\Gamma_i$ for $1 \leq i \leq p$. Then each $|\Gamma_i|$ is a tree that has $w_{E_b}$ as an endpoint. For $1 \leq i \leq p$, let $A_i = \cup \{ cl(F) : F \in \mathcal{F}$ and $w_F$ is a vertex of $\Gamma_i \}$. Then $A = A_1 \cup A_2 \cup \ldots \cup A_p$ and $A_i \cap A_j \subset E_b$ for $i \neq j$. Since $v \notin E_b$ and $(A_i - E_b) \cap (A_j - E_b) = \emptyset$ for $i \neq j$, then $v$ belongs to at most one of the sets $A_i$. Since $p \geq 3$, we can assume, after reindexing, that $v \notin A_1 \cup A_2$. Note that $E_b \subset cl(F_i) \subset A_i$ for $i = 1, 2$.

Now for $i = 1, 2$, let $\Psi_i$ be a subcomplex of $\Gamma_i$ with the following four properties.

- $w_{E_b}$ is an endpoint of $\Psi_i$, and $\Psi_i$ contains the vertex $w_{F_i}$ and the edge joining $w_{E_b}$ to $w_{F_i}$.

- For every vertex of $\Psi_i$ of the form $w_E$ where $E \in \mathcal{E} - \{E_b\}$, $\deg(w_E, \Psi_i) = 2$.

- For every vertex of $\Psi_i$ of the form $w_F$ where $F \in \mathcal{F}$, $\Psi_i$ contains every edge of $\Gamma$ that emanates from $w_F$.

- $|\Psi_i|$ is connected.

To construct $\Psi_i$, orient the edges of $\Gamma_i$ away from $w_{E_b}$, include the edge from $w_{E_b}$ to $w_{F_i}$ and all the edges emanating from $w_{F_i}$, but only one outgoing edge emanating from a vertex $w_E$ that is at the terminal end of an edge emanating from $w_{F_i}$. Continue in this fashion to construct $\Psi_i$ inductively: when one encounters a new vertex at the terminal end of an edge, add only one outgoing edge if the vertex is of the form $w_E$ ($E \in \mathcal{E}$), but add all the outgoing edges if the vertex is of the form $w_F$ ($F \in \mathcal{F}$). This process must terminate because $\Gamma_i$ is a finite complex. It can't stop at a vertex of the form $w_E$ ($E \in \mathcal{E}$) because all these vertices have degree $\geq 2$ in $\Gamma$, since $Q = \Delta \cup_f J$ has no exposed edges. The process only terminates when it has adjoined all the vertices it can reach of degree $\geq 2$, and only has access to vertices of $\Gamma$ of the form $w_F$ ($F \in \mathcal{E}$) that are endpoints of $\Gamma$. The process then adjoins these endpoints to $\Psi_i$ until they exhausted.

For $i = 1, 2$, let $R_i = \cup \{ cl(F) : F \in \mathcal{F}$ and $w_F$ is a vertex of $\Psi_i \}$. Then $F_i \subset R_i \subset A_i$. Moreover, each $R_i$ is a union of finitely many disks such that the intersection of any two of the disks is either empty or an arc lying in the boundary of each, and the intersection of any three distinct disks is empty. Thus, each $R_i$ is a compact 2-manifold. The strong deformation retraction of $A$ onto $\Gamma$ restricts to a strong deformation retraction of $A_i$ onto $\Gamma_i$ which, in turn, restricts to a strong deformation retraction of $R_i$ onto $\Psi_i$. Thus, each $R_i$ is contractible and, hence, is a disk. Clearly, $E_b$ is an arc in $\partial R_i$.

Observe that $fr(A) = \cup \{ cl_Q(\partial F) : F \in \mathcal{F} \}$. Also observe that for $i = 1, 2$, $\partial R_i = ( \cup \{ cl_Q(\partial F) : F \in \mathcal{F}$ and $w_F$ is a vertex of $\Psi_i \} ) \cup in(E_b)$. Thus, $\partial R_i - in(E_b) \subset fr(A)$.

Let $R = R_1 \cup R_2$. Since $R_1$ and $R_2$ are disks and $R_1 \cap R_2 = \partial R_1 \cap \partial R_2 = E_b$, then $R$ is a disk, $\partial R = (\partial R_1 - in(E_b)) \cup (\partial R_2 - in(E_b))$ and $in(E_b) \subset in(R)$. Hence, $\partial R \subset fr(A)$. Since $Q = A \cup B$ and $B$ is a closed subset of $Q$, then $fr(A) \subset B$. Therefore, $\partial R \subset B$. Furthermore, $v \notin R$.

Since $b \in in(E_b) \subset in(R)$, we can choose an arc $H$ in $in(E_b)$ such that $b \in in(H)$, and for $i = 1, 2$, we can choose a disk $G_i \subset cl(F_i)$ such that $G_i \cap \partial(cl(F_i)) = H$. Hence, $H \subset \partial G_1 \cap \partial G_2$ and $G_1 \cap G_2 = H$. Thus, $G_i - H \subset F_i$. Therefore, $G = G_1 \cup G_2$ is a disk in $in(R)$ such that $v \notin G$ and $b \in in(H) \subset in(G)$. Let $j : \partial G \to Q - \{b\}$ denote the inclusion map. Since $R - in(G)$ is an annulus in $Q - \{b\}$ with boundary components $\partial R$ and $\partial G$, then $j$ is homotopic in $Q - \{b\}$ to a map $\hat{k} : \partial G \to Q - \{b\}$ such that $\hat{k}(\partial G) \subset \partial R$. Hence, the induced maps $j_* : H_1(\partial G) \to H_1(Q - \{b\})$ and $\hat{k}_* : H_1(\partial G) \to H_1(Q - \{b\})$ are equal. Since $\partial R \subset B \subset Q - \{b\}$, then $\hat{k}_*$ factors through the finite group $H_1(B)$. Since $H_1(\partial G) \approx \mathbb{Z}$, it follows that $\hat{k}_*$ and, hence, $j_*$ are not monomorphisms, contradicting Lemma 4. This concludes the proof of Subcase **2i)**.

Recall that $E_b$ is the element of $\mathcal{E}$ that contains $b$. Subcase **2i)** dealt with the situation in which $v \notin E_b$. We now assume $v \in E_b$. Then $w_{E_b} = v$. We will break this situation into two subcases. First we identify $St(v,Q)$ and $Lk(v,Q)$ with subsets of $Q$ with more specificity than previously. Recall that in this situation, $E_b = \cup_{i=1}^{n} E_i$ where $E_i = E_b \cap J_i \in \mathcal{E}_i$ and $v \in in(E_i)$ for $1 \leq i \leq n$, and $b \in in(E_1)$. Furthermore, $St(v,J) = St(v,Q) \cap J$ and $Lk(v,J) = Lk(v,Q) \cap J$ are identified with $\cup_{i=1}^{n} E_i = E_b$ and $\cup_{i=1}^{n} \partial E_i$, respectively. These identifications make $St(v,Q) \cap J_1 = E_1$ and $Lk(v,Q) \cap J_1 = \partial E_1$. The components of $q^{-1}(E_b)$ form a finite pairwise disjoint collection of arcs $\{ L_1, L_2, \ldots , L_k \}$ in $\partial \Delta$ such that $f(\cup_{j=1}^{k} \partial L_j) = \cup_{i=1}^{n} \partial E_i = Lk(v,J)$. Furthermore, for $1 \leq j \leq k$, there is a point $u_j \in in(L_j)$ such that $q^{-1}(\{v\}) \cap L_j = \{u_j\}$ and such that $q$ maps each component of $L_j - \{u_j\}$ homeomorphically onto one of the components of $E_b - \{v\}$. (Either $q$ embeds $L_j$ in $E_b$ so that its image is an arc joining two distinct points of $\cup_{i=1}^{n} \partial E_i$, or $q$ "folds" $L_j$ at $u_j$ and maps each component of $L_j - \{u_j\}$ homeomorphically onto the same component of $E_1 - \{v\}$.) As we saw at the beginning of Subcase **2*i***), there are k distinct elements $F_1$, $F_2$, $\ldots$ , $F_k$ of $\mathcal{F}$ such that $L_j \subset \partial P(F_j)$ for $1 \leq j \leq k$. (Then an edge of $\Gamma$ joins $w_{E_b} = v$ to $w_{F_j}$ for $1 \leq j \leq k$.) For $1 \leq j \leq k$, let $L_j'$ be an arc in the disk $P(F_j)$ such that $\partial L_j' = \partial L_j$ and $in(L_j') \subset in(P(F_j))$ for $1 \leq j \leq k$. Then our description near the beginning of Section 4 of how to identify $Lk(v,Q)$ with a subset of Q implies that $Lk(v,Q)$ can be identified with $q(\cup_{j=1}^{k} L_j')$. Observe that for $1 \leq j \leq k$: either $q$ embeds $L_j$ in $E_b$ in which case $q(L_j')$ is an arc in $Lk(v,Q)$ joining two distinct points of $Lk(v,J) = \cup_{i=1}^{n} \partial E_i$, or $q$ "folds" $L_j$ and $q(L_j')$ is a simple closed curve in $Lk(v,Q)$ that contains a single point of $Lk(v,J) = \cup_{i=1}^{n} \partial E_i$. Recall that our choice at the beginning of the proof of Case 2 guarantees that each of the two points of $\partial E_1 = Lk(v,Q) \cap J_1$ lie in simple closed curves in $Lk(v,Q)$ and that $f^{-1}(J_1 - \{v\})$ has at least 3 components. Since $b \in in(E_1)$ and $b \neq v$, then there is an endpoint $e$ of $E_1$ such that $b$ and $e$ lie in the same component of $E_1 - \{v\}$. (Observe, that the component of $E_1 - \{v\}$ that contains b and has endpoint e is also the component of $E_b - \{v\}$ that contains v and has endpoint e.) Since $e \in \partial E_1$, then $e$ lies in a simple closed curve in $Lk(v,Q)$. We can now formulate Subcases **2*ii***) and **2*iii***).

**Subcase 2*ii*).** For some j between 1 and k, $q(L_j')$ is a simple closed curve that contains e.

**Subcase 2*iii*).** There is no j between 1 and k such that $q(L_j')$ is a simple closed curve that contains e.

Assume we are in Subcase **2*ii*)**. By reindexing $\{ L_1', L_2', \ldots , L_k' \}$, we can assume $q(L_1')$ is a simple closed curve in $Lk(v,Q)$ that contains e. Then $u_1 \in in(L_1)$, $q^{-1}(\{v\}) \cap L_1 = \{u_1\}$, $q$ maps each of the two components of $L_1 - \{u_1\}$ homeomorphically onto the component of $E_1 - \{v\}$ that contains e and b. $F_1$ is the element of $\mathcal{F}$ such that $L_1 \subset \partial P(F_1)$; moreover, $L_1' \subset P(F_1)$, $\partial L_1' = \partial L_1$ and $in(L_1') \subset in(P(F_1))$. Furthermore, q maps $P(F_1) - L_1$ homeomorphically onto $cl(F_1) - E_1$ and q *zips up* $L_1$ so that $q(L_1)$ is an arc in the disk $cl(F_1)$, $\partial q(L_1) = \{ e, v \}$, $e \in \partial(cl(F_1))$ and $q(L_1) - \{e\} \subset in(cl(F_1))$.

Recall that in this situation, $w_{E_b} = v$. Let $\Gamma_1$ be the subgraph of $\Gamma$ such that $|\Gamma_1|$ is the closure of the component of $|\Gamma| - \{w_{E_b}\}$ that contains $w_{F_1}$. Then $|\Gamma_1|$ is a tree and $w_{E_b}$ is an endpoint of $\Gamma_1$. Let $A_1 = \cup \{ cl(F) : F \in \mathcal{F}$ and $w_F$ is a vertex of $\Gamma_1 \}$.

As in Subcase **2i)**, let $\Psi_1$ be a subcomplex of $\Gamma_1$ with the following four properties.

- $w_{E_b}$ is an endpoint of $\Psi_1$, and $\Psi_1$ contains the vertex $w_{F_1}$ and the edge joining $w_{E_b}$ to $w_{F_1}$.

- For every vertex of $\Psi_1$ of the form $w_E$ where $E \in \mathcal{E} - \{E_b\}$, $\deg(w_E, \Psi_1) = 2$.

- For every vertex of $\Psi_1$ of the form $w_F$ where $F \in \mathcal{F}$, $\Psi_1$ contains every edge of $\Gamma$ that emanates from $w_F$.

- $|\Psi_1|$ is connected.

Let $R_1 = \cup \{ cl(F) : F \in \mathcal{F}$ and $w_F$ is a vertex of $\Psi_1 \}$. Then $F_1 \subset R_1 \subset A_1$. Moreover, $R_1$ is a union of finitely many disks such that the intersection of any two of the disks is either empty or an arc lying in the boundary of each, and the intersection of any three distinct disks is empty. Thus, $R_1$ is a compact 2-manifold. The strong deformation retraction of $A$ onto $\Gamma$ restricts to a strong deformation retraction of $A_1$ onto $\Gamma_1$ which, in turn, restricts to a strong deformation retraction of $R_1$ onto $\Psi_1$. Thus, $R_i$ is contractible and, hence, is a disk.

Recall that $fr(A) = \cup \{ cl_Q(\partial F) : F \in \mathcal{F} \}$. Observe that since $q(L_1) - \{e\} \subset in(cl(F_1))$, then $\partial R_1 = ( \cup \{ cl_Q(\partial F) : F \in \mathcal{F}$ and $w_F$ is a vertex of $\Psi_1 \} )$. Hence, as in Subcase **2i)**, $\partial R_1 \subset fr(A) \subset B$.

Since $b \in in(q(L_1)) \subset in(F_1)$, we can choose an arc H in $in(q(L_1))$ such that $b \in in(H)$, and we can choose disks $G_1$ and $G_2$ in $in(cl(F_1))$ such that $G_i \cap q(L_1) = H \subset \partial G_i$ for $i = 1, 2$ and $G_1 \cap G_2 = H$. Therefore, $G = G_1 \cup G_2$ is a disk in $in(cl(F_1)) \subset in(R_1)$ such that $v \notin G$ and $b \in in(H) \subset in(G)$. Let $j : \partial G \to Q - \{b\}$ denote the inclusion map. Since $R_1 - in(G)$ is an annulus in $Q - \{b\}$ with boundary components $\partial R_1$ and $\partial G$, then $j$ is homotopic in $Q - \{b\}$ to a map $\hat{k} : \partial G \to Q - \{b\}$ such that $\hat{k}(\partial G) \subset \partial R_1$. Hence, the induced maps $j_* : H_1(\partial G) \to H_1(Q - \{b\})$ and $\hat{k}_* : H_1(\partial G) \to H_1(Q - \{b\})$ are equal. Since $\partial R_1 \subset B \subset Q - \{b\}$, then $\hat{k}_*$ factors through the finite group $H_1(B)$. Since $H_1(\partial G) \approx \mathbb{Z}$, it follows that $\hat{k}_*$ and, hence, $j_*$ are not monomorphisms, contradicting Lemma 4. This concludes the proof of the Theorem in Subcase **2ii)**.

Finally assume we are in Subcase **2iii**). By hypothesis, there is a simple closed curve in Lk(v,Q) that contains e, but no individual $q(L_j')$ is such a simple closed curve for $1 \leq j \leq k$. Therefore, there is an r between 2 and k such that, after reindexing, we can assume $q(L_1')$, $q(L_2')$, ... , $q(L_r')$ are arcs in Lk(v,Q) whose union is a simple closed curve that contains e. Moreover we can index the points of $\cup_{j=1}^r q(\partial L_j')$ as $c_1$, $c_2$, ... , $c_r$ = $c_0$ so that $q(\partial L_j') = \{ c_{j-1}, c_j \}$ for $1 \leq j \leq r$ and e = $c_0$ = $c_r$. Recall that for $1 \leq j \leq r$, $L_j$ is an arc in $\partial\Delta$ which is a component of $q^{-1}(E_b)$, $F_j \in \mathcal{F}$ such that $L_j$ is a component of $P(F_j) \cap \partial\Delta$, $L_j'$ is an arc in $P(F_j)$ such that $\partial L_j = \partial L_j'$ and $in(L_j') \subset in(P(F_j))$. (Let $F_0 = F_r$.) Recall that $F_1$, $F_2$, ... , $F_r$ must be distinct, because if $F_j = F_{j'}$ for $j \neq j'$, then there are 2 distinct edges of $\Gamma$ joining $w_{F_j} = w_{F_{j'}}$ to $w_{E_b}$, creating a loop in the tree $|\Gamma|$. Recall that $E_b$ is a cone with vertex v and base $\cup_{i=1}^n \partial E_i$, and that each $q(L_j)$ is the (unique) arc in $E_b$ joining two points of $\cup_{i=1}^n \partial E_i$. Thus, $\{ c_1, c_2, ... , c_r \} \subset \cup_{i=1}^n \partial E_i$; and if we let $M_j$ denote the unique arc in $E_b$ joining v to the point $c_j$ for $0 \leq j \leq r$, then $q(L_j) = M_{j-1} \cup M_j$. Thus, $L_j$ is a component of $P(F_j) \cap \partial\Delta$ that q maps into $E_b$. Since b lies in the component of $E_b - \{v\}$ with end point e = $c_0$ = $c_r$ and b $\neq$ e, then b $\in in(M_0) = in(M_r)$. If $L'$ is a component of $P(F_j) \cap \partial\Delta$ such that $L' \neq L_j$, then $q(L') \cap E_b = \emptyset$, because otherwise (as above) there would be 2 distinct edges of $\Gamma$ joining $w_{F_j}$ to $w_{E_b}$, creating a loop in the tree $|\Gamma|$. Thus, $cl(F_j) \cap E_b = q(P(F)) \cap E_b = q(L_j) = M_{j-1} \cup M_j$ for $1 \leq j \leq r$. Also for $j \neq j'$, $cl(F_j)$ and $cl(F_{j'})$ can't both intersect an element E of $\mathcal{E}$ such that $E \neq E_b$, because that would give rise to an edge loop in $\Gamma$ with vertices $w_{E_b}$, $w_{F_j}$, $w_E$, $w_{F_{j'}}$, creating a loop in the tree $|\Gamma|$. Hence, $cl(F_j) \cap cl(F_{j'}) = (cl(F_j) \cap E_b) \cap (cl(F_{j'}) \cap E_b) = (M_{j-1} \cup M_j) \cap (M_{j'-1} \cup M_{j'})$. Thus, $cl(F_j) \cap cl(F_{j'}) = M_{\min\{j,j'\}}$ if $j' = j \pm 1$ and $cl(F_j) \cap cl(F_{j'}) = \{v\}$ if either $j' \leq j - 2$ of $j' \geq j + 2$. Finally, note that since each $q(L_j')$ is an arc rather than a simple closed curve, then q embeds each $L_j$ in $E_b \subset J$ and, consequently, $q(L_j)$ is an arc in $\partial cl(F_j)$. Thus, $M_{j-1}$ and $M_j$ are arcs in $\partial cl(F_j)$.

Recall that $w_{E_b} = v$. Let $1 \leq j \leq r$. Let $\Gamma_j$ be the subcomplex of $\Gamma$ such that $|\Gamma_j|$ is the closure of the component of $\Gamma - \{ w_{E_b} \}$ that contains $w_{F_j}$. Then $|\Gamma_j|$ is a tree, $w_{E_b}$ is an endpoint of $\Gamma_j$, and $\Gamma_j$ contains an edge that joins $w_{E_b}$ to $w_{F_j}$. Let $\mathcal{F}_j = \{ F \in \mathcal{F} : w_F$ is a vertex of $\Gamma_j \}$ and let $A_j = \cup \{ cl(F) : F \in \mathcal{F}_j \}$.

For $1 \leq j \leq r$, as in Subcase **2i**), let $\Psi_j$ be a subcomplex of $\Gamma_j$ with the following four properties.

- $w_{E_b}$ is an endpoint of $\Psi_j$, and $\Psi_j$ contains the vertex $w_{F_j}$ and the edge joining $w_{E_b}$ to $w_{F_j}$.

- For every vertex of $\Psi_j$ of the form $w_E$ where $E \in \mathcal{E} - \{E_b\}$, $\deg(w_E, \Psi_j) = 2$.

- For every vertex of $\Psi_j$ of the form $w_F$ where $F \in \mathcal{F}$, $\Psi_j$ contains every edge of $\Gamma$ that emanates from $w_F$.

- $|\Psi_j|$ is connected.

For $1 \le j \le r$, let $R_j = \cup \{ cl(F) : F \in \mathcal{F}$ and $w_F$ is a vertex of $\Psi_j \}$. Then $cl(F_j) \subset R_j \subset A_j$. Moreover, $R_j$ is a union of finitely many disks such that the intersection of any two of the disks is either empty or an arc lying in the boundary of each, and the intersection of any three distinct disks is empty. Thus, $R_j$ is a compact 2-manifold. The strong deformation retraction of $A$ onto $\Gamma$ restricts to a strong deformation retraction of $A_j$ onto $\Gamma_j$ which, in turn, restricts to a strong deformation retraction of $R_j$ onto $\Psi_j$. Thus, $R_j$ is contractible and, hence, is a disk.

For $1 \le t \le r$, let $S_t = \cup_{j=1}^{t} R_j$. We claim:

- $S_t$ is a disk for $1 \le t \le r$,

- $\cup_{j=1}^{t-1} in(M_j) \subset in(S_t)$ for $2 \le t \le r-1$ and $(\cup_{j=1}^{r} in(M_j)) \cup \{v\} \subset in(S_r)$, and

- $M_0 \cup M_t \subset \partial S_t$ for $1 \le t \le r-1$, $\partial S_t = (\cup_{j=1}^{t} \partial R_j) - (\cup_{j=1}^{t-1} in(M_j))$ for $2 \le t \le r-1$ and $\partial S_r = \cup_{j=1}^{r} (\partial R_j - in(M_{j-1} \cup M_j))$.

To prove this claim we begin with some basic observations. For $E \in \mathcal{E}$ and $F \in \mathcal{F}$, $cl(F) \cap E \neq \emptyset$ if and only if $w_E$ and $w_F$ are joined by an edge of $\Gamma$. Also for $F$ and $F' \in \mathcal{F}$, $cl(F) \cap cl(F') \neq \emptyset$ if and only if there is an $E \in \mathcal{E}$ such that $w_E$ is joined to $w_F$ and $w_{F'}$ by two edges of $\Gamma$. Let $R_j' = \cup \{ cl(F) : F \in \mathcal{F} - \{F_j\}$ and $w_F$ is a vertex of $\Psi_j \}$. Clearly, no edge of $\Gamma$ joins $w_{E_b}$ to a vertex $w_F$ of $\Psi_j$ where $F \in \mathcal{F} - \{F_j\}$. Hence, $E_b \cup R_j' = \emptyset$. Since $cl(F_j)$ and $R_j = cl(F_j) \cup R_j'$ are disks, $cl(F_j) \cap E_b = M_{j-1} \cup M_j \subset \partial cl(F_j)$, and $R_j'$ is a closed subset of $R_j$ that is disjoint from $E_b$, then $R_j \cap E_b = M_{j-1} \cup M_j \subset \partial R_j$. Hence, $R_j \cap in(M_{j'}) = \emptyset$ for $j' \le j-2$ and $j' \ge j+1$. If $j \neq j'$, $F \in \mathcal{F} - \{F_j\}$ such that $w_F$ is a vertex of $\Psi_j$, and $F' \in \mathcal{F}$ such that $w_{F'}$ is a vertex of $\Psi_{j'}$, then there is no $E \in \mathcal{E}$ such that edges of $\Gamma$ join $w_E$ to $w_F$ and to $w_{F'}$. Hence, $cl(F) \cap cl(F') = \emptyset$. It follows that $R_j' \cap R_{j'} = \emptyset$. Thus, $R_j \cap R_{j'} = cl(F_j) \cap cl(F_{j'})$ for $j \neq j'$. Consequently, $R_j \cap R_{j'} = M_{min\{j,j'\}}$ if $j' = j \pm 1$ and $R_j \cap R_{j'} = \{v\}$ if either $j' \le j-2$ of $j' \ge j+2$.

We now turn to the proof of the previous claim. The proof is inductive. Clearly $S_1 = R_1$ is a disk and $M_0 \cup M_1 \subset \partial R_1 = \partial S_1$. Let $1 \le t \le r-2$ and assume $S_t = \cup_{j=1}^{t} R_j$ satisfies the assertions of the claim. Then $S_t \cap R_{t+1} = (\cup_{j=1}^{t} R_j) \cap R_{t+1} = \cup_{j=1}^{t} (R_j \cap R_{t+1}) = \{v\} \cup M_t = M_t$ where $M_t$ is an arc lying in $\partial S_t$ and in $\partial R_{t+1}$. Hence, $S_{t+1} = S_t \cup R_{t+1}$ is a disk such that $in(S_{t+1}) = in(S_t) \cup in(M_t) \cup in(R_{t+1})$ and $\partial S_{t+1} = (\partial S_t - in(M_t)) \cup (\partial R_{t+1} - in(M_t))$. Since $\cup_{j=1}^{t-1} in(M_j) \subset in(S_t)$ by inductive hypothesis, then clearly $\cup_{j=1}^{t} in(M_j) \subset in(S_{t+1})$. Since $M_0 \subset \partial S_t$ by inductive hypothesis, $M_0 \cap in(M_t) = \emptyset$, $M_{t+1} \subset \partial R_{t+1}$ and $M_{t+1} \cap in(M_t) = \emptyset$, then $M_0 \cup M_{t+1} \subset \partial S_t$. Since $R_{t+1} \cap in(M_j) = \emptyset$ for $j \le t-1$, then $\partial R_{t+1} - in(M_t) = \partial R_{t+1} - (\cup_{j=1}^{t} in(M_j))$. By inductive hypothesis, $\partial S_t = (\cup_{j=1}^{t} \partial R_j) - (\cup_{j=1}^{t-1} in(M_j))$. Hence, $(\partial S_t - in(M_t)) = (\cup_{j=1}^{t} \partial R_j) - (\cup_{j=1}^{t} in(M_j))$. Therefore,

$$\partial S_{t+1} = (\partial S_t - in(M_t)) \cup (\partial R_{t+1} - in(M_t)) =$$

$$((\cup_{j=1}^{t} \partial R_j) - (\cup_{j=1}^{t} in(M_j))) \cup (\partial R_{t+1} - (\cup_{j=1}^{t} in(M_j))) = \cup_{j=1}^{t+1} \partial R_j) - (\cup_{j=1}^{t} in(M_j)).$$

This establishes the claim for $1 \leq t \leq r - 1$. Finally consider $S_r = S_{r-1} \cup R_r$. $S_{r-1} \cap R_r = (\cup_{j=1}^{r-1} R_j) \cap R_r = \cup_{j=1}^{r-1} (R_j \cap R_r) = M_0 \cup \{v\} \cup M_{r-1} = M_0 \cup M_{r-1}$ where $M_0 \cup M_{r-1} = M_{r-1} \cup M_r$ is an arc lying in $\partial S_{r-1}$ and in $\partial R_r$. Hence, $S_r = S_{r-1} \cup R_r$ is a disk such that $in(S_r) = in(S_{r-1}) \cup in(M_0 \cup M_{r-1}) \cup in(R_r)$ and $\partial S_r = (\partial S_{r-1} - in(M_0 \cup M_{r-1})) \cup (\partial R_r - in(M_0 \cup M_{r-1}))$. Since $\cup_{j=1}^{r-2} in(M_j) \subset in(S_{r-1})$ by inductive hypothesis and $in(M_{r-1}) \cup in(M_r) \cup \{v\} = in(M_0) \cup in(M_{r-1}) \cup \{v\} = in(M_0 \cup M_{r-1})$, then clearly $(\cup_{j=1}^{r} in(M_j)) \cup \{v\} \subset in(S_r)$. Since $R_r \cap in(M_j) = \emptyset$ for $j \leq r - 2$ and $in(M_0 \cup M_{r-1}) = in(M_{r-1}) \cup in(M_r) \cup \{v\}$, then $\partial R_r - in(M_0 \cup M_{r-1}) = \partial R_r - ((\cup_{j=1}^{r} in(M_j)) \cup \{v\})$. By inductive hypothesis, $\partial S_{r-1} = (\cup_{j=1}^{r-1} \partial R_j) - (\cup_{j=1}^{r-2} in(M_j))$. Hence, $\partial S_{r-1} - in(M_0 \cup M_{r-1}) = ((\cup_{j=1}^{r-1} \partial R_j) - (\cup_{j=1}^{r-2} in(M_j))) - (in(M_{r-1}) \cup in(M_r) \cup \{v\}) = (\cup_{j=1}^{r-1} \partial R_j) - ((\cup_{j=1}^{r} in(M_j)) \cup \{v\})$. Therefore,

$$\partial S_r = (\partial S_{r-1} - in(M_0 \cup M_{r-1})) \cup (\partial R_r - in(M_0 \cup M_{r-1})) =$$

$$((\cup_{j=1}^{r-1} \partial R_j) - ((\cup_{j=1}^{r} in(M_j)) \cup \{v\})) \cup (\partial R_r - ((\cup_{j=1}^{r} in(M_j)) \cup \{v\})) =$$

$$(\cup_{j=1}^{r} \partial R_j) - ((\cup_{j=1}^{r} in(M_j)) \cup \{v\}) = \cup_{j=1}^{r} (\partial R_j - ((\cup_{j'=1}^{r} in(M_{j'})) \cup \{v\})).$$

Since $R_j \cap in(M_{j'}) = \emptyset$ for $j' \leq j - 2$ and $j' \geq j + 1$ and $in(M_{j-1}) \cup in(M_j) \cup \{v\} = in(M_{j-1} \cup M_j)$, then $\partial R_j - ((\cup_{j'=1}^{r} in(M_{j'})) \cup \{v\}) = \partial R_j - in(M_{j-1} \cup M_j)$. Consequently,

$$\partial S_r = \cup_{j=1}^{r} (\partial R_j - in(M_{j-1} \cup M_j)).$$

This completes the proof of the claim.

As we observed previously, $fr(A) = \cup \{ cl_Q(\partial F) : F \in \mathcal{F} \}$. Also observe that for $1 \leq j \leq r$, $\partial R_j = ( \cup \{ cl_Q(\partial F) : F \in \mathcal{F}$ and $w_F$ is a vertex of $\Psi_j \} ) \cup in(M_{j-1} \cup M_j)$. Thus, $\partial R_j - in(M_{j-1} \cup M_j) \subset fr(A)$. The previous claim implies $\partial S_r = \cup_{j=1}^{r} (\partial R_j - in(M_{j-1} \cup M_j))$. Thus, $\partial S_r \subset fr(A)$. Hence, $\partial S_r \subset B$.

Recall that $in(M_r) \subset J - \{v\}$. Since $b \in in(M_r) \subset in(S_r) - \{v\}$, we can choose an arc $H$ in $in(M_r)$ such that $b \in in(H)$, and we can choose disks $G_1$ and $G_2$ in $in(S_r)$ such that $G_i \cap J = H \subset \partial G_i$ for $i = 1, 2$ and $G_1 \cap G_2 = H$. Therefore, $G = G_1 \cup G_2$ is a disk in $in(S_r)$ such that $v \notin G$ and $b \in in(H) \subset in(G)$. Let $j : \partial G \to Q - \{b\}$ denote the inclusion map. Since $S_r - in(G)$ is an annulus in $Q - \{b\}$ with boundary components $\partial S_r$ and $\partial G$, then $j$ is homotopic in $Q - \{b\}$ to a map $\hat{k} : \partial G \to Q - \{b\}$ such that $\hat{k}(\partial G) \subset \partial S_r$. Hence, the induced maps $j_\star : H_1(\partial G) \to H_1(Q - \{b\})$ and $\hat{k}_\star : H_1(\partial G) \to H_1(Q - \{b\})$ are equal. Since $\partial S_r \subset B \subset Q - \{b\}$, then $\hat{k}_\star$ factors through the finite group $H_1(B)$. Since $H_1(\partial G) \approx \mathbb{Z}$, it follows that $\hat{k}_*$ and, hence, $j_*$ are not monomorphisms, contradicting Lemma 4. This concludes the proof of the Theorem in Subcase **2*iii*).** ∎